%%% TXP %%% BEGIN prytz.tex %%%%%%%%%%%%%%%%%%%%%%%%%%%%%%%%%%%%%%%%%%%%%%%
% Geometry of the Prytz Planimeter
% by Robert L. Foote
% Wabash College
% footer@wabash.edu
% 11 Aug 98.
% To appear in Reports on Mathematical Physics.

% 23 pages.  12 figures are contained in two *.gif files at 600 dpi.

% Compile with AmS-TeX.
\input amstex
\input epsf
\documentstyle{amsppt}
\magnification\magstep1
\NoRunningHeads
\TagsOnRight
%\nopagenumbers

% Input macro files.
%%% TXP %%% BEGIN prytz.mac %%%%%%%%%%%%%%%%%%%%%%%%%%%%%%%%%%%%%%%%%%%%%%%
% Macros for Prytz planimeter paper.  These are all used.
\vsize8.9truein\hsize6.5truein

\define \ie {{\it i.e.}}
\define \eg {{\it e.g.}}

\redefine \l{\ell}
\define\sech{\operatorname{sech}}
\define\csch{\operatorname{csch}}
\define\st{|}
\define \iso  {\approx}                 % Isomorphism relation.
\define\comp{\mathbin{\ssize\circ}}     % Function composition.

%\define\qand{\quad\text{and}\quad}
\define\qqand{\qquad\text{and}\qquad}

% Blackboard bold.
\define \R{{\Bbb R}}                   % Set of real numbers.
\define \C{{\Bbb C}}                   % Set of complex numbers.

\define \cross{\times}                  % Cartesian product or cross product.
\define \bdy {\partial}                 % Boundary.
\redefine \phi {\varphi}
\define \union    {\cup}
\define \To#1{\overset{#1}\to\longrightarrow}
           % For putting the name of the function over the arrow bewteen the
           % domain and range, e.g.,  $A \To f B$.
           % Note: The construction @>#1>> puts #1 too high.
\define \inv  {^{-1}}    % For the inverse of a function or matrix.

\define\smallpmatrix{\(\smallmatrix}
\define\endsmallpmatrix{\endsmallmatrix\)}

% Nested parentheses.
\define \({\left(}
\define \){\right)_{\vphantom{\leavevmode\lower1pt\hbox{x}}}}

\define \lowsub#1{_{\vphantom{\tilde{#1}}#1}}
                           % Sometimes subscripts are not low enough, eg,
                           % \mu_L for the left-invariant Mauer-Cartan form.
                           % Using \mu\lowsub L makes the subscript taller
                           % by an invisible tilde, pushing it lower.

%\input /pctex/macros/Calculus.mac
\define \der#1#2{\frac{d#1}{d#2}}        % Derivative of #1 wrt #2.
\define \pder#1#2{\frac{\partial#1}{\partial#2}} % Parital derivative.
\define \derat#1#2#3{\left.\frac{d#1}{d#2}\right|_{#3}}
                         % Derivative of #1 wrt #2 evaluated at #3.
\define \pderat#1#2#3{\left.\frac{\partial#1}{\partial#2}\right|_{#3}}
                         % Parital derivative of #1 wrt #2 at #3.

%\input /pctex/macros/CxAna.mac
\define  \conj   {\overline}
\redefine \Re    {\operatorname{Re}}  % As in the real part of something.
\redefine \Im    {\operatorname{Im}}  % As in the imaginary part of something.

% Additional macros
\define\a{\alpha}
\redefine\b{\beta}
\redefine\d{\delta}
\define\g{\gamma}
\define\G{\Gamma}
\define\gt{\tilde\gamma}
\define\s{\sigma}
\define\sT{\sigma\lowsub{\!T}}
\redefine\O{\Omega}
\define\bdyO{{\partial\Omega}}
\redefine\o{\omega}
\define\ot{\tilde\omega}
\define\M{\Cal M_0(S^1)}  % Group of Mobius transformations acting on S^1.
%\define\m{\frak m(S^1)}   % Lie algebra of \M.
\define\SU{SU(1,1)}
\define\su{\frak{su}(1,1)} 
\define\Diff{\operatorname{Diff}(S^1)}
\define\X{\frak X(S^1)}
\define\z{\zeta}
\define\tensor{\otimes}
\define\dirsum{\oplus}
\define\Ad{\operatorname{Ad}}
\define\tr{\operatorname{tr}}
%%% TXP %%% END prytz.mac %%%%%%%%%%%%%%%%%%%%%%%%%%%%%%%%%%%%%%%%%%%%%%%%%

\topmatter
 \title Geometry of the Prytz Planimeter \endtitle
 \author  Robert L. Foote  \endauthor
 \affil   Wabash College  \endaffil
 \address Department of Mathematics \& Computer Science,
          Crawfordsville, IN 47933 \endaddress
 \email footer\@wabash.edu 
        \newline\indent
        {\it Web page:\/} http://persweb.wabash.edu/facstaff/footer/ \endemail
 \thanks
   The author would like to acknowledge helpful discussions on this topic with 
   Felix Albrecht, Larry Bates, Richard Bishop, Lance Drager, and Jeff Lee.
 \endthanks
 \thanks
   To appear in {\sl Reports on Mathematical Physics.}
 \endthanks
 \subjclass 70F25,   % Nonholonomic systems
            53B15,   % Local Geometry, Other connections
            53A99,   % Classical DG, None of the above but in this section
            53C65    % Global DG, Integral geometry 
 \endsubjclass
 \abstract 
   The Prytz planimeter is a simple example of a system governed by a non-holonomic 
   constraint. It is unique among planimeters in that it measures something more 
   subtle than area, combining the area, centroid and other moments of the region 
   being measured, with weights depending on the length of the planimeter. As a 
   tool for measuring area, it is most accurate for regions that are small relative 
   to its length.

   The configuration space of the planimeter is a non-principal circle bundle acted 
   on by $\SU\,(\iso SL(2,\R))$. The motion of the planimeter is realized as parallel 
   translation for a connection on this bundle and for a connection on a principal 
   $\SU$-bundle. The holonomy group is $\SU$. As a consequence, the planimeter is an 
   example of a system with a phase shift on the circle that is not a simple rotation.

   There is a qualitative difference in the holonomy when tracing large regions as 
   opposed to small ones. Generic elements of $\SU$ act on $S^1$ with two fixed 
   points or with no fixed points. When tracing small regions, the holonomy acts 
   without fixed points. Menzin's conjecture states (roughly) that if a planimeter of 
   length $\l$ traces the boundary of a region with area $A > \pi\l^2$, then it 
   exhibits an asymptotic behavior and the holonomy acts with two fixed 
   points, one attracting and one repelling. This is obvious if the region is a disk, 
   and intuitively plausible if the region is convex and $A\gg\pi\l^2$. A proof of 
   this conjecture is given for a special case, and the conjecture is shown to imply 
   the isoperimetric inequality.
 \endabstract
\endtopmatter

\document
\baselineskip14pt plus 1pt minus 1pt
\overfullrule0pt

%\baselineskip=24pt % Controls spacing between lines. 12pt is default.

%%% TXP %%% BEGIN Prytz123.tex %%%%%%%%%%%%%%%%%%%%%%%%%%%%%%%%%%%%%%%%%%%%

% Geometry of the Prytz Planimeter, Intro and Sections 1,2,3.

A planimeter is a mechanical instrument used to determine the area of a region in the plane.
The user moves a ``tracing point'' around the boundary of the region. When the tracing point
returns to its starting point, some feature of the instrument does not quite return to its 
initial position. This can be interpreted as holonomy, and is related to the area of the region.

%\vskip2.5truein  %picture of polar and linear planimeters
%\vfill

%\newpage

Most planimeters consist of a rod, one end of which is the tracing point $T$. A wheel is attached
to the rod that partially rolls and partially slides on the paper as the tracer point is moved.  
The most familiar such instrument is the polar planimeter (Figure~1(a)), invented by Jacob Amsler
in 1854, in which the end of the rod opposite
the tracer point is hinged to a secondary rod, restricting its motion to a circle.
In the linear planimeter (Figure~1(b)) the other end of the tracing rod is restricted to move 
along a linear track.
The ``roll'' of the wheel is recorded on a scale, which reads out the area when the tracer point 
returns to its starting point. These work because the wheel mechanically integrates a 1-form that
differs from
$\frac12(-y\,dx + x\,dy)$ by an exact form (see \cite{Fo1} for details and other references).
\bigskip
\centerline{\epsfxsize.9\hsize\epsfbox{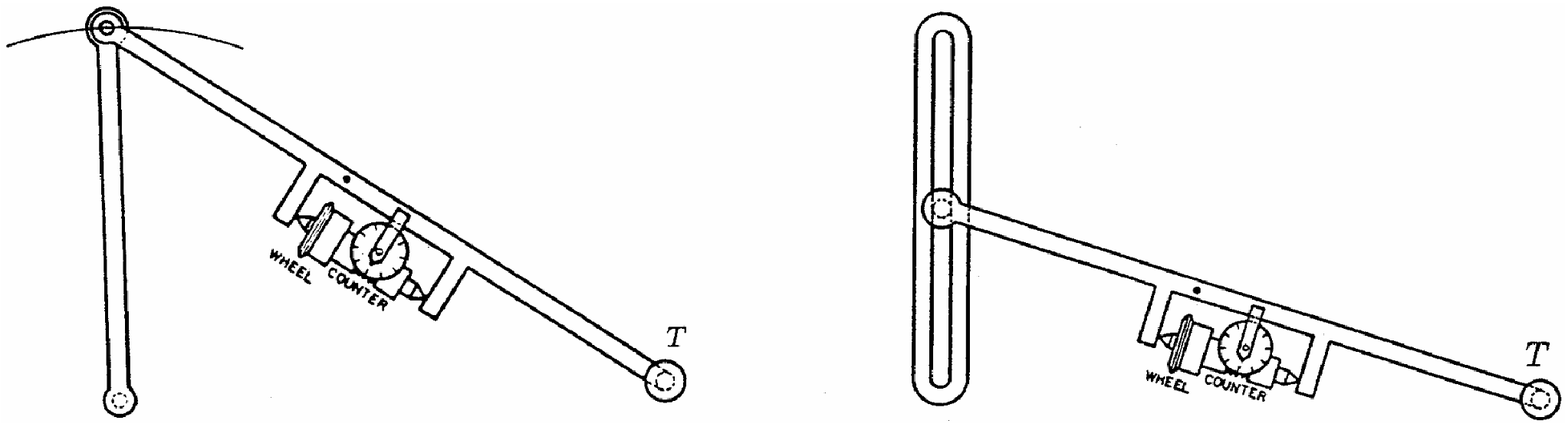}}
%\vskip6truecm  %picture of polar and linear planimeters
%\bigskip
\centerline{Figure 1(a). Polar Planimeter. \hfil\hfil Figure 1(b). Linear Planimeter.}
\bigskip

The Prytz planimeter (Figure~2), in contrast, contains no internally moving parts. It 
consists simply of
a rod with its ends bent at right angles. One end, the tracer point $T$, is sharpened to a 
point. The other end, $C$, is sharpened to a chisel edge parallel to the rod. The chisel edge is 
usually slightly rounded, making it look similar to a hatchet, and consequently the device is 
also known as a ``hatchet planimeter.'' It was invented in about 1875 by Holger Prytz, a Danish 
cavalry officer and mathematician, as an economical and simple alternative to Amsler's
planimeter. Prytz referred to it as a ``stang planimeter,'' ``stang'' being Danish for ``rod.'' 
An amusing account of the history of the Prytz planimeter is given by Pedersen \cite{Pe}. For a 
very complete history of planimeters through 1894, see Henrici \cite{He}.

%\vfill
%\vskip7truecm %picture of Prytz planimeter and standard tractrix w/ area shaded
\centerline{\epsfxsize.45\hsize\epsfbox{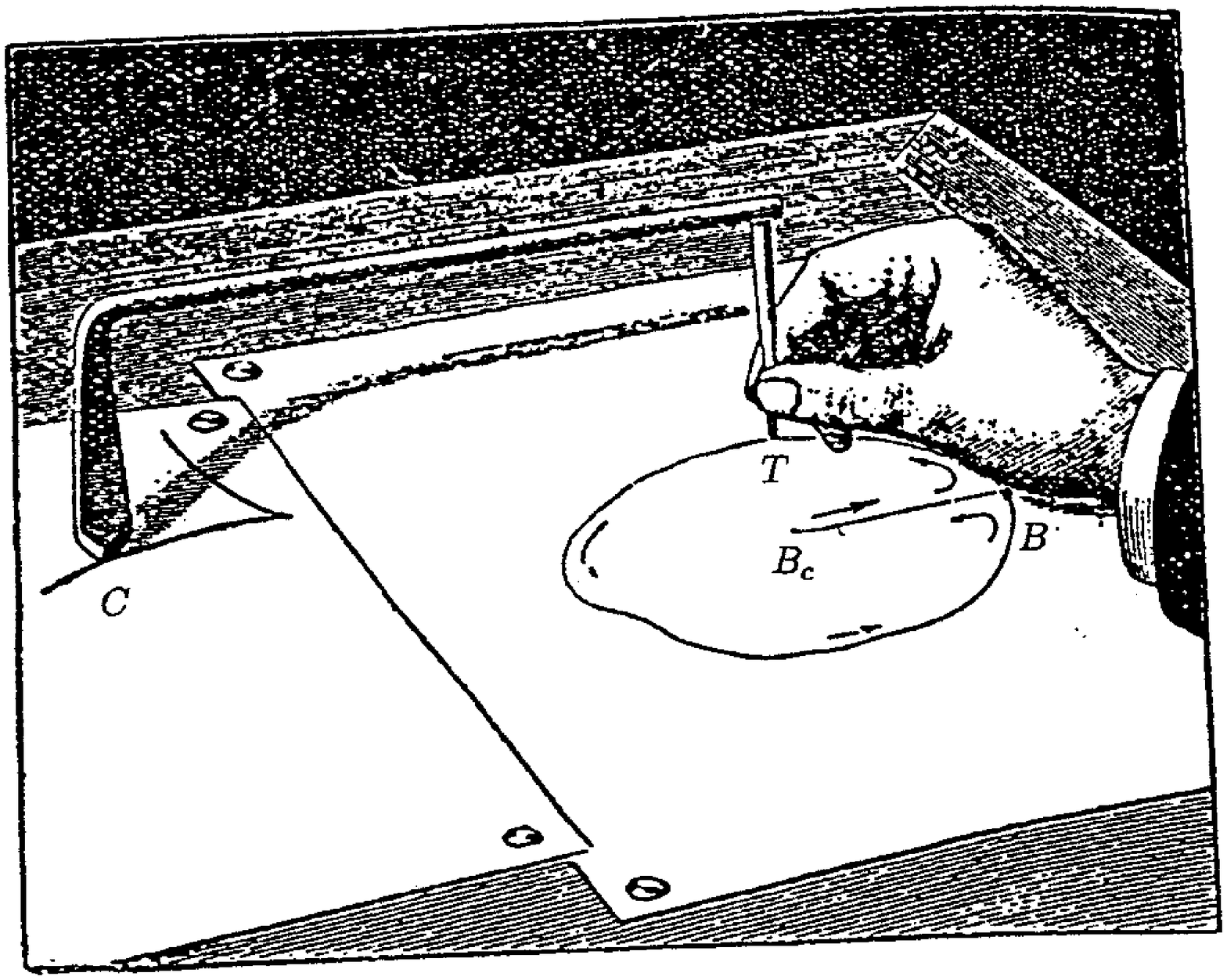}\qquad\epsfxsize.45\hsize\epsfbox{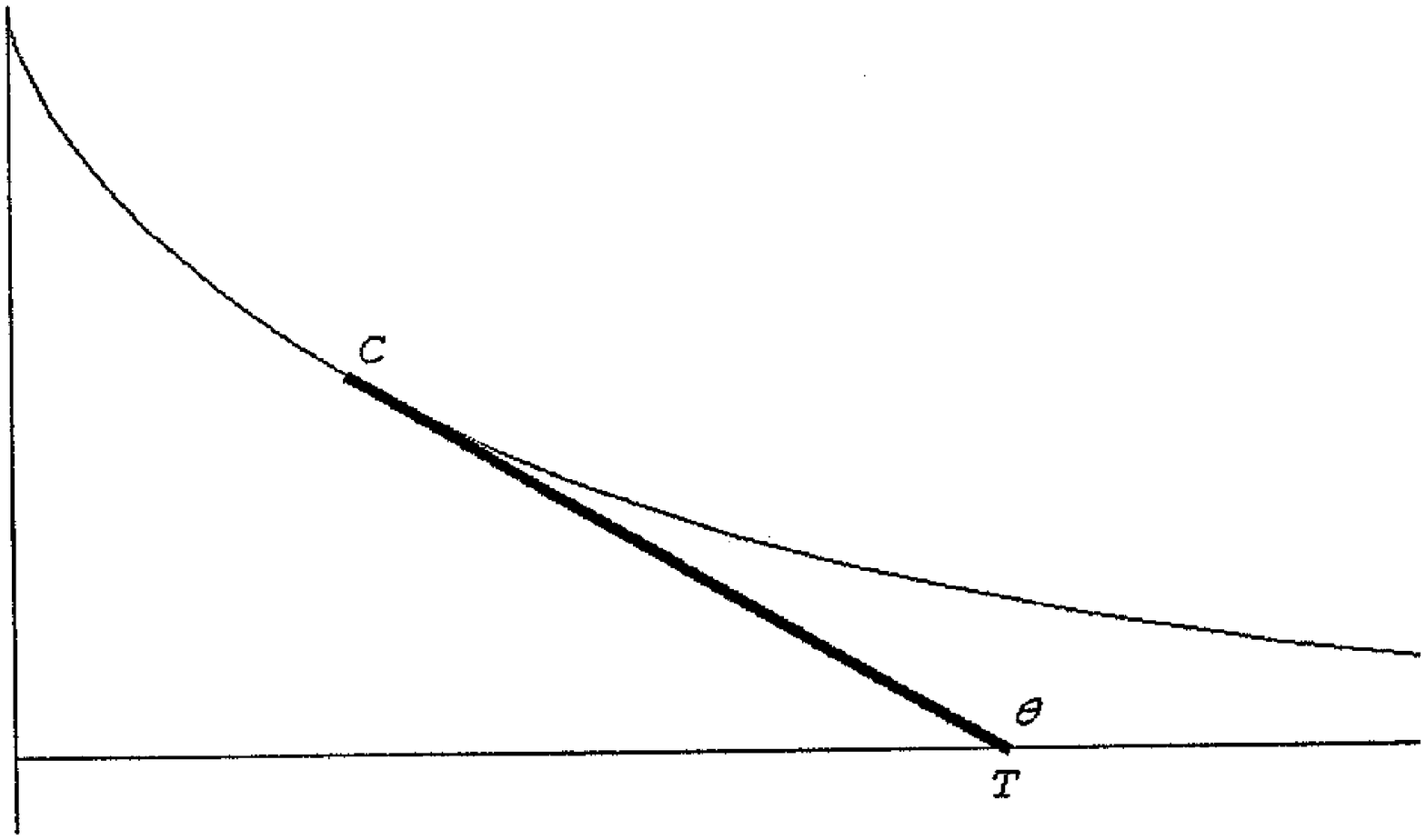}}
%\bigskip
\centerline{Figure 2. Prytz planimeter. \hfil\hfil Figure 3. Standard tractrix.}
\bigskip

To use the planimeter, grasp it at the end with the tracer point, keeping its ends 
perpendicular to the plane. Move the tracer point along some curve taking care not to 
apply any torque to it. As $T$ moves along the curve, $C$ follows a pursuit curve that
is always tangent to the rod. The paths of $T$ and $C$ are similar to those followed 
by the front and rear wheels of a bicycle. Note that if $T$ is moved along a straight 
line, the path of $C$ is a standard tractrix (Figure~3). Thus when $T$ moves along 
an arbitrary curve $\gamma$, the path of $C$ is called a tractrix of $\gamma$.  
Tractrices of various curves, particularly of circles, are considered in 
\cite{L,Mo,Pou}.

%\vfillneg

%\newpage

%{\rightskip3truein  % picture of zig-zag path
To measure the area of a region $\O$, start $T$ at some base point $B$ on $\bdyO$, and
note the initial position of $C$.  As $T$ moves around $\bdyO$, $C$ describes a 
zig-zag path (Figure~4), which Richard Bishop has likened to the motion of parking a 
car.  When $T$ gets back to $B$, the chisel edge $C$ does not return to its starting 
point---it has undergone a displacement. Multiply the displacement $\s$ of $C$ by the 
length $\l$ of the rod and you have the area. Actually it's only an approximation, but
the error turns out to have a very nice geometric description.%\par}

We have implicitly assumed here that $\bdyO$ consists of a simple, closed curve, so
$\O$ is bounded and simply connected. To precisely describe the behavior of the
planimeter, we will also need to assume that $\bdyO$ is piecewise $C^1$. These
assumptions about $\O$ will be made for the rest of the paper without further comment.

%\vfill\vfill\vfill
%\vskip7.5truecm
%\bigskip
\epsfxsize.9\hsize
\centerline{\epsfbox{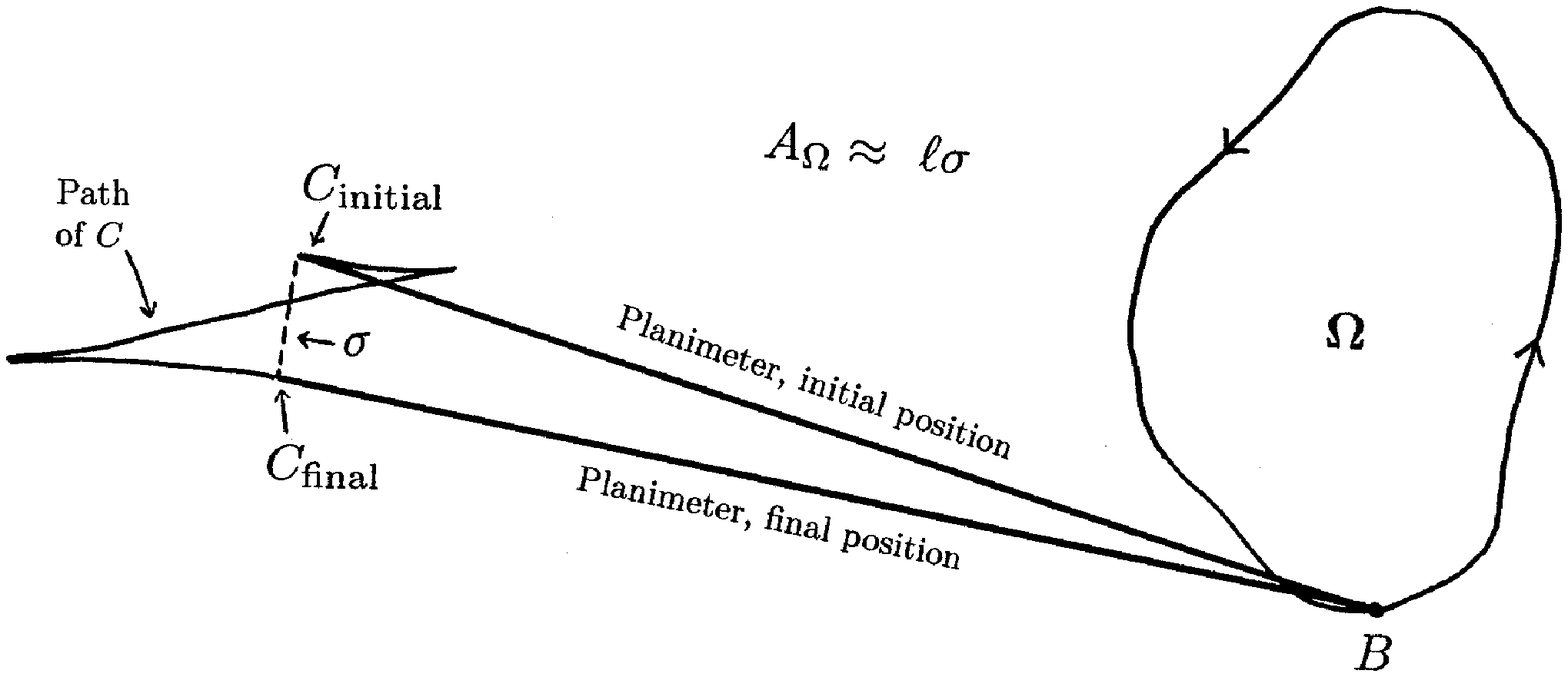}}
\centerline{Figure 4. Motion of the Prytz planimeter.}
\bigskip

Despite its simplicity and inherent inaccuracy as a measuring tool---in fact, {\sl because\/} of 
its inherent inaccuracy---the Prytz planimeter is more interesting mathematically than the polar
planimeter.

After giving the elementary theory of the Prytz planimeter, we will discuss Menzin's 
conjecture on how the planimeter behaves when tracing large regions and its connection
to the isoperimetric inequality, the analytic approach of the inventor and one of his 
contemporaries, the behavior of the planimeter as a non-holonomic system---most 
specifically as parallel translation for a connection on a principal $\SU$-bundle and
for the induced connection on an associated circle bundle---and
finally, a proof of a special case of Menzin's conjecture.

\heading 1. Elementary Theory \endheading

Let $p$ and $q$ be distinct points in $\R^2$, and consider the segment
joining them.  Let $N$ be the unit vector perpendicular to the segment so that $q-p$ and $N$
form a positively-oriented frame.

%\vskip3.5truecm
%\bigskip
\centerline{\epsfxsize.45\hsize\epsfbox{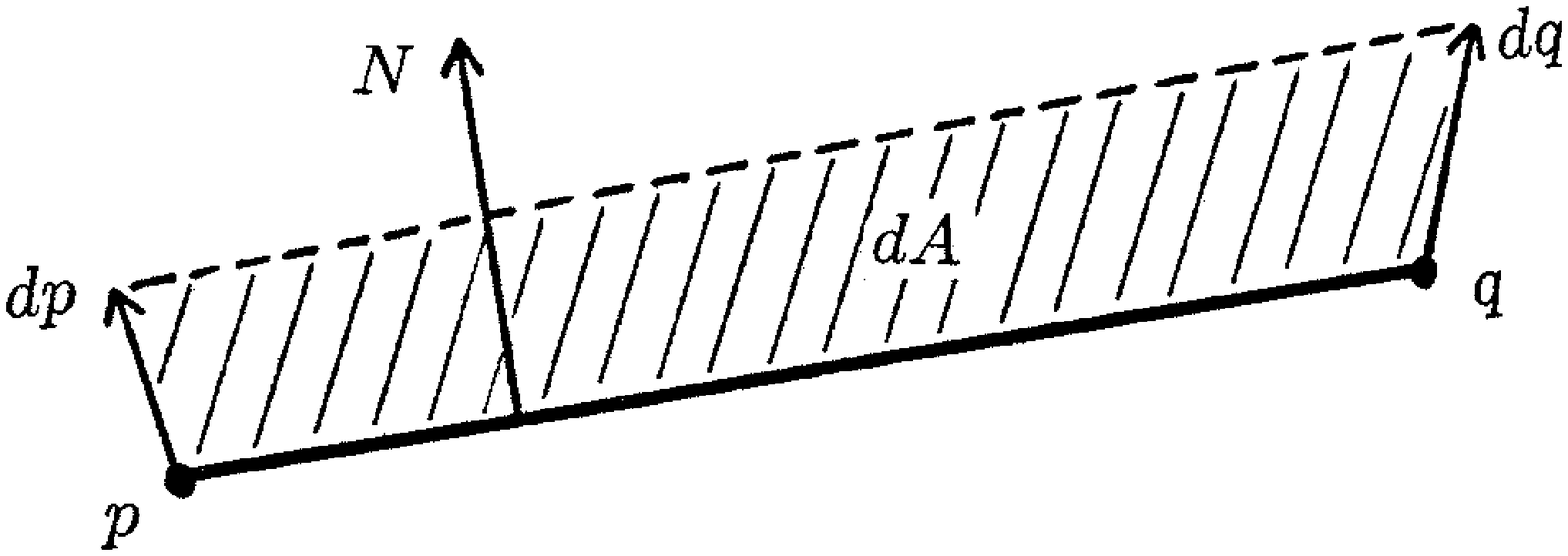}\qquad\epsfxsize.45\hsize\epsfbox{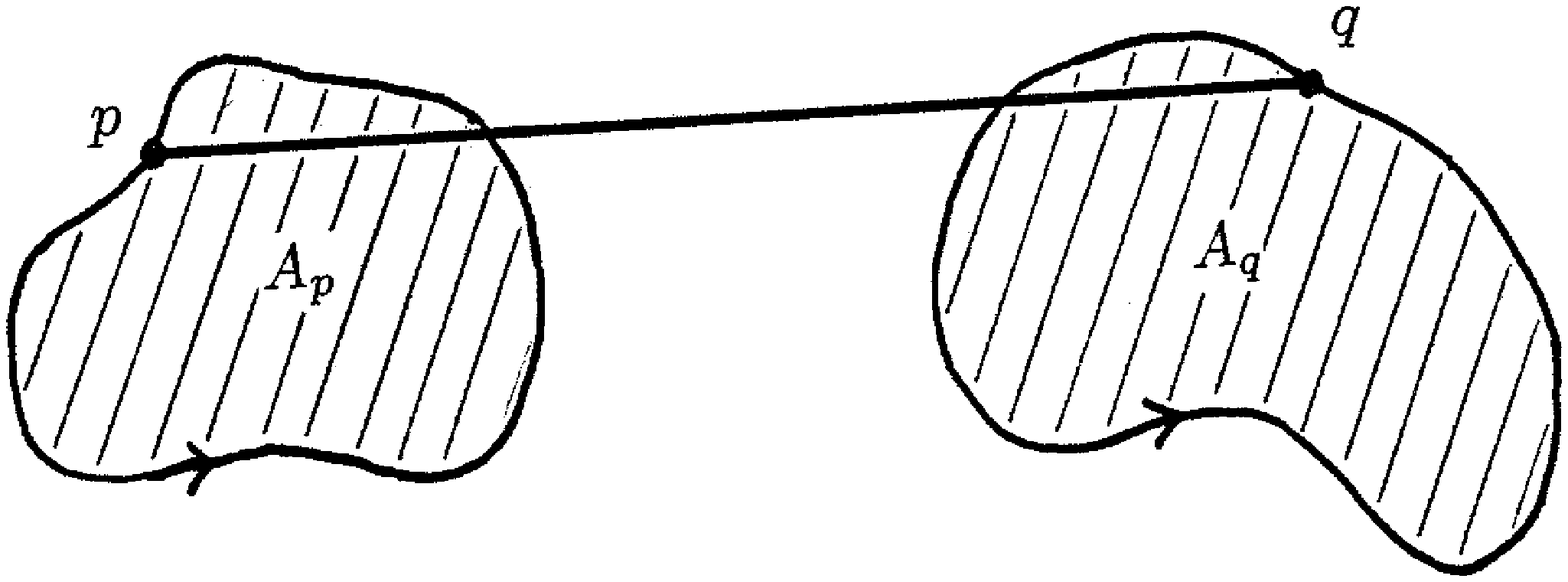}}
\centerline{Figure 5. \hfil\hfil Figure 6.}
\bigskip
%\vskip2truein  %picture of inf movement and of moving segment thm
%\vfill

%\newpage

If the segment moves slightly, an infinitesimal oriented area is swept out (Figure~5).
If $m = \frac12(p+q)$ is the midpoint of the segment, then this area is given by 
$dA = \l N{\cdot}dm = \frac12\l N{\cdot}(dp + dq)$, where $\l$ is the length of the 
segment (which can be variable), and $dp$, $dq$, and $dm$ are the infinitesimal
displacements of $p$, $q$, and $m$, respectively. Note that $N$ points in the 
direction of motion for which the area is counted positively. 

Let $d\s = N {\cdot} dp$ be the component of $dp$ in the direction of $N$. Noting that
the infinitesimal rotation of the segment is $d\theta = \frac1\l N{\cdot}(dq-dp)$, 
where $\theta$ is the angle of inclination of the segment relative to some fixed 
direction, the expression for the infinitesimal area can be written as 
$dA = \l\,d\s + \frac12 \l^2\,d\theta$.

These expressions for $dA$ are intuitively plausible. If $dp = dq$, then 
$dA = \l N{\cdot}dm = \frac12\l N{\cdot}(dp + dq)$ is just the area of a 
parallelogram.  If $dp$ and $dq$ are parallel to the segment, they are orthogonal to
$N$ and no area is swept out.  If the segment rotates about its midpoint, then 
$dm = 0$ and the oriented areas swept out by the two halves of the segment cancel.  
These motions account for all four dimensions of the configuration space of the moving
segment, and so any infinitesimal motion of the segment is a linear combination of 
these.  A similar intuitive argument can be made for 
$dA = \l\,d\s + \frac12 \l^2\,d\theta$. More rigorously, the expression 
$dA = \l N{\cdot}dm$ is the integrand in the two-dimensional version of Guldin's 
formula for the area swept out by a moving segment \cite{Cou,Fo1}. 

The formula for infinitesimal rotation (which will be needed later) is easily obtained
by differentiating $N{\cdot}(p-q) = 0$ and noting that $dN/d\theta$ is the unit vector
$u=\frac1\l(p-q)$.  We get 
$0 = u\,d\theta{\cdot}(p-q) + N{\cdot}(dp-dq) = \l\,d\theta + N{\cdot}(dp-dq)$.

%\vfillneg

%\newpage

If the endpoints of the moving segment trace out curves, the expressions for 
$dA$ can be integrated, yielding the total oriented area swept out. In 
particular, if the endpoints trace out closed curves (Figure~6), it is easy to show 
\cite{Fo1} that
the oriented area swept out is $A_q - A_p$, where $A_p$ and $A_q$ are the 
oriented areas enclosed by curves traced out by $p$ and $q$, respectively.

The moving segment becomes a Prytz planimeter if we fix its length $\l$ and use the 
endpoint $q$ as the tracer point $T$.  The behavior of the planimeter is then given 
simply by the non-holonomic restriction $d\s = 0$, that is, the infinitesimal movement
of the chisel point $C$ has no component in the direction of $N$. As a consequence, 
$dA = \frac12 \l^2\,d\theta$ and the total oriented area swept out by the planimeter 
at any moment is $A = \frac12 \l^2\,\Delta\theta$, where $\Delta\theta$ is the net 
rotation of the planimeter.  This gives a simple proof that the area under the 
standard tractrix (Figure~3) is $\frac14\pi\l^2$.

Now consider what happens when the tracer point is moved around the boundary of a region $\O$ 
in the positive direction, starting and ending at some base point $B \in \bdy\O$ 
(Figure~4). The chisel 
point $C$ does not return to 
its original position. For regions that aren't too large relative to $\l$ (and for which the
planimeter gives the most accurate results, as we will see), the chisel point comes to rest 
at a point close to its original position.  Now imagine the user fixing the tracer point $T$ at 
$B$ and rotating the planimeter about $B$ (this motion violates the non-holonomic constraint)
so that $C$ moves along a circle back to its initial position 
(we will call this circle---with center $B$ and radius $\l$---the ``initial circle''). As both 
endpoints would now have
traced out closed curves, the oriented area swept out is $A_q - A_p = A_\O - A_\g$, where $A_\O$ 
is the area of $\O$ and $A_\g$ is the oriented area enclosed by the curve $\g$ formed by the
zig-zag 
path of $C$ during the normal use of the planimeter {\sl plus\/} the arc of the initial circle 
from the final position of $C$ back to its initial position.  The oriented area 
$A_q - A_p$ also equals the integral of
$dA = \l\,d\s + \frac12 \l^2\,d\theta$.  As the initial and final angles of the planimeter are the
same, $d\theta$ integrates to 0.  Along the part of $\g$ followed by the chisel we have
$d\s = 0$, but along the part of $\g$ that is the arc of the initial circle,
$d\s$ integrates to $\s$, the length of that arc. Thus we have
$$
A_\O - A_\g = \l\s, \qquad\text{or}\qquad A_\O = \l\s + A_\g,
$$
and we see that the error made by the approximation $A_\O \approx \l\s$ is $A_\g$. In a typical 
use
of the planimeter, the curve $\g$ encloses a number of ``triangular'' regions. These regions are 
generally small compared to $\O$. Moreover the boundaries of those regions inside the initial
circle are traced with the opposite orientation of those outside, and so their oriented areas
have opposite signs in their sum $A_\g$. Keeping this in mind, a good starting position is one
that will result in the chisel edge spending part of its time inside the initial circle and
part of its time outside. For example, one could start with the planimeter perpendicular to a 
line that roughly bisects the region.

Another way to minimize the error, suggested by several authors 
\cite{Ba,Cr,Hi,K,Pr1-4,Sa,St}, is
to draw a line segment from $B$ to the centroid $B_c$ of the
region.  Instead of starting the tracer point at $B$, start at $B_c$.
Trace along the segment $BB_c$ to $B$, then around the curve, and finally
back along $BB_c$, stopping at $B_c$, as suggested in Figure~2.  This, of course, leads to the
problem of locating the centroid, which is at least as hard as computing the
area! In practice, one simply guesses. Prytz and Hill are the only authors who give any
argument that the centroid should be used as the base point (outlined in Section~3). 
The other authors simply appeal to
this as a geometrically plausible way to balance the triangular regions inside and outside the 
initial circle.

Similar arguments can be made to explain how other planimeters work, and have appeared 
in a variety of forms dating from the early 1800's (see \cite{Fo1} for references). 
The argument given above to explain the behavior of the Prytz planimeter in particular
is due to Henrici \cite{He}, who was the first to give a common theoretical setting 
for most of the planimeters invented up to 1894. This argument also appears in 
a paper by Kriloff \cite{K}, and may have been discovered independently by him.

It is interesting that this simple geometric argument can be made without referring to 
the exact relationship between the movement of the tracer point and the angle $\theta$
between the planimeter rod and some fixed direction. This relationship is the focus of
much of the rest of the paper. Recall the expression for the infinitesimal rotation of
a moving line segment, $\l\,d\theta = N{\cdot}(dq - dp)$, which for the planimeter 
becomes $\l\,d\theta = N{\cdot}dq$, since $d\s = N{\cdot}dp = 0$. If the tracer point 
$T=q$ follows the curve $\g\:\R\to\R^2$, the differential equation governing $\theta$ 
is then $\l\der\theta t = N{\cdot}\g'(t)$. If $\theta$ measures the angle between the
planimeter and the $x$-axis, \ie, 
$p-q = \l(\cos\theta\,\bold i + \sin\theta\,\bold j)$,
the equation becomes $\l\der\theta t = \sin\theta \der xt - \cos\theta \der yt$.
If the tracer point moves along the $x$-axis with $x=t$, as in Figure~3, we get
$\l\der\theta x = \sin\theta$, the solution of which is easily seen to be
$$
\tan\frac\theta 2 = Ae^{x/\l}.  \tag1
$$
This can be rewritten as \,$\tan\theta = -\csch(a + x/\l)$, where $a = \log A$.
If the planimeter starts perpendicular to the $x$-axis when $x=0$, then $a=0$, and
so \,$\tan\theta = -\csch(x/\l)$. This is the slope of the line tangent to
the tractrix when that curve is parameterized by the point where the tangent meets
the axis (see, \eg, \cite{Cox}):
$$
\tau(x) = \big(x - \l\tanh(x/\l),\ \l \sech(x/\l)\big).  \tag2
$$

Variations of the Prytz planimeter have been made, the most notable ones being by Goodman and
Scott. As we have seen, the length of the arc along the initial circle between the initial and 
final locations of the chisel edge is the important quantity to measure, and not simply the 
distance between these points, although the latter is usually adequate. The variations
due to Goodman and Scott allow the 
direct measurement of this arc. Goodman \cite{G} incorporated a curved scale with radius $\l$ 
into the planimeter rod, so the scale lies along the initial circle when applied to the points
(Figure~7).
Scott's variation \cite{Sc} implements the idea in the discussion above 
that brings the chisel end of the
planimeter back to its initial position. He put a wheel next to the chisel edge with its axis 
parallel to the rod. This wheel does not contact the paper during the tracing of $\bdyO$, but
rides slightly above it. When 
the tracing point returns to the base point, the planimeter is tilted slightly, bringing the 
wheel into contact with the paper. The user then rotates the planimeter around the tracer point.
The wheel rolls along the initial circle between the points, measuring the arc-length.
Both inventors, particularly Scott, seemed to believe that they were addressing the 
cause of the instrument's error, but the proof above and the analysis by Prytz and 
Hill in Section~3 show that the error is more complicated (and more interesting) 
than this.

%\newpage
%\ \vskip4truecm  % picture of Goodman's planimeter
\bigskip
\centerline{\epsfxsize.9\hsize\epsfbox{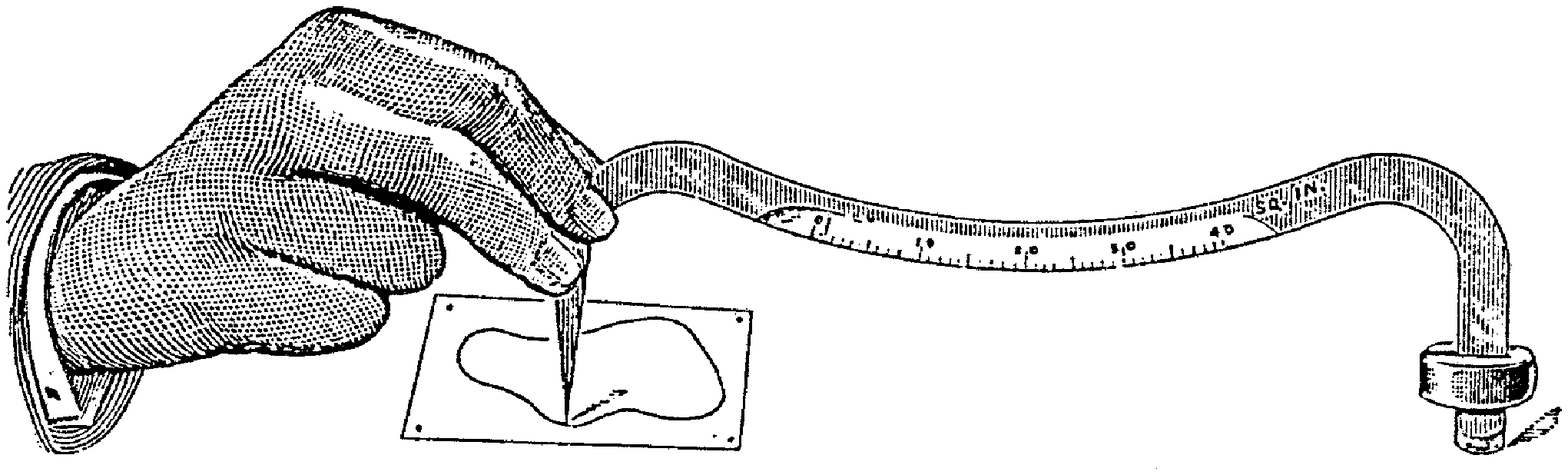}}
\centerline{Figure 7. Goodman's planimeter.}
\bigskip

%\newpage

The planimeters of Prytz, Goodman, and Scott were all marketed (the original Prytz 
planimeter was produced by the firm of Knudsen in Copenhagen \cite{Pe}). The 
additional details of the latter two, however, defeated Prytz's original purpose of 
economy relative to Amsler's planimeter. In reaction to Goodman and Scott, Prytz 
\cite{Pr3} advised engineers ``rather than use the `improved stang planimeters,' let a
country blacksmith make them a copy of the original instrument.'' A collector has sent
the author two photographs, one of a Prytz planimeter manufactured in the Netherlands 
and one of a Goodman planimeter, which are posted on the author's web page (URL at the
end of this paper). The author would be very interested to learn of other instruments 
of this type that still exist.

\heading 2. Menzin's Conjecture and the Isoperimetric Inequality \endheading

When the planimeter traces the boundary of a polygonal region $\O$ with edges that are
large compared to $\l$, the length of the planimeter, it is intuitively clear that the
path of $C$ asymptotically approaches a particular tractrix that is closed, and that 
when the planimeter follows this closed tractrix, it makes a complete rotation. This 
also happens when $\O$ is a disk of radius $R > \l$; it is easily shown that the 
closed tractrix is a circle of radius $\sqrt{R^2 - \l^2}$. Furthermore, it is
intuitively plausible for large convex regions. After much experimentation, Menzin 
\cite{Me} conjectured that this happens whenever $A_\O > \pi\l^2$. The last section 
contains a proof of Menzin's conjecture when $\bdyO$ is a parallelogram.

In spite of its title (having to do with Amsler's planimeter), the paper by Morley 
\cite{Mo} goes into considerable detail on the tractrices of circles. In addition to
the case where $R>\l$, he notes that when $R = \l$ the asymptotic tractrix reduces to 
the center of the circle, but it is attractive only from one 
side. When $R < \l$ each tractrix is made up of regularly spaced cusps.

It turns out that Menzin's conjecture implies the isoperimetric inequality.  To see this, note
that when $C$ follows the closed tractrix, $dA = \l\,d\s + \frac12\l^2\,d\theta$ integrates to 
$\pi\l^2$, since 
$d\s = 0$ and the planimeter makes a complete rotation.  Since $C$ traces a closed path, the
total area swept out is
$A_\O - A_C$, where $A_C$ is the area enclosed by the path of $C$ (which does {\sl not\/} 
include an arc of the initial circle).  Furthermore, note that $dA$ can be written as
$dA = \l\,d\sT - \frac12\l^2\,d\theta$, where $d\sT = N {\cdot} dq$ is the component of $dq$
(the infinitesimal displacement of the tracer point) in the direction of $N$. This 
expression for $dA$ 
integrates to $\l\sT - \pi\l^2$, where $\sT = \int d\sT$. (Note that if a wheel is mounted
on the rod at $T$ similar to the wheel on a polar or linear planimeter, then $\sT$ is the total
signed distance the wheel rolls.) Thus we have
$$
\pi\l^2 = A_\O - A_C = \l\sT - \pi\l^2.
$$
By Menzin's conjecture this holds under the assumption that $A_\O > \pi\l^2$.  In the limiting 
case where $A_\O = \pi\l^2$, we have $A_C = 0$. A little algebra yields
$\sT^2 = 4\pi A_\O$. Noting that $d\sT$ measures only a component of the infinitesimal
distance the tracer point moves ($d\sT \le ds$), we have $\sT \le L_{\bdyO}$, where 
$L_{\bdyO}$ is the length of $\bdyO$.  The isoperimetric inequality follows:
$$
L_{\bdyO}^2 \ge \sT^2 = 4\pi A_\O.
$$
Furthermore, $L_{\bdyO}^2 = 4\pi A_\O$ implies $d\sT = ds$, that is, the tracer point only moves 
in the direction perpendicular to the rod. In this case the chisel edge does not move at all, 
and the tracer point describes a circle of radius $\l$. Thus the isoperimetric inequality is an
equality only if $\bdy\O$ is a circle.

For other connections between planimeters and isoperimetric inequalities, see 
\cite{Fo1\&2}.

\heading 3. Analysis by Prytz and Hill \endheading

%{\rightskip2.5truein  %picture from Hill. I redrew it.
Prytz and Hill studied the ``stang'' planimeter analytically, rather than geometrically.
Prytz's original paper \cite{Pr1} was published anonymously in Danish under the pseudonym `Z'
(see \cite{Pe} for the story behind this). His subsequent short notes on this are very tersely 
written \cite{Pr2\&4}.  Hill's account is much more readable \cite{Hi}.

%\vskip6truecm
%\bigskip
\centerline{\epsfxsize.5\hsize\epsfbox{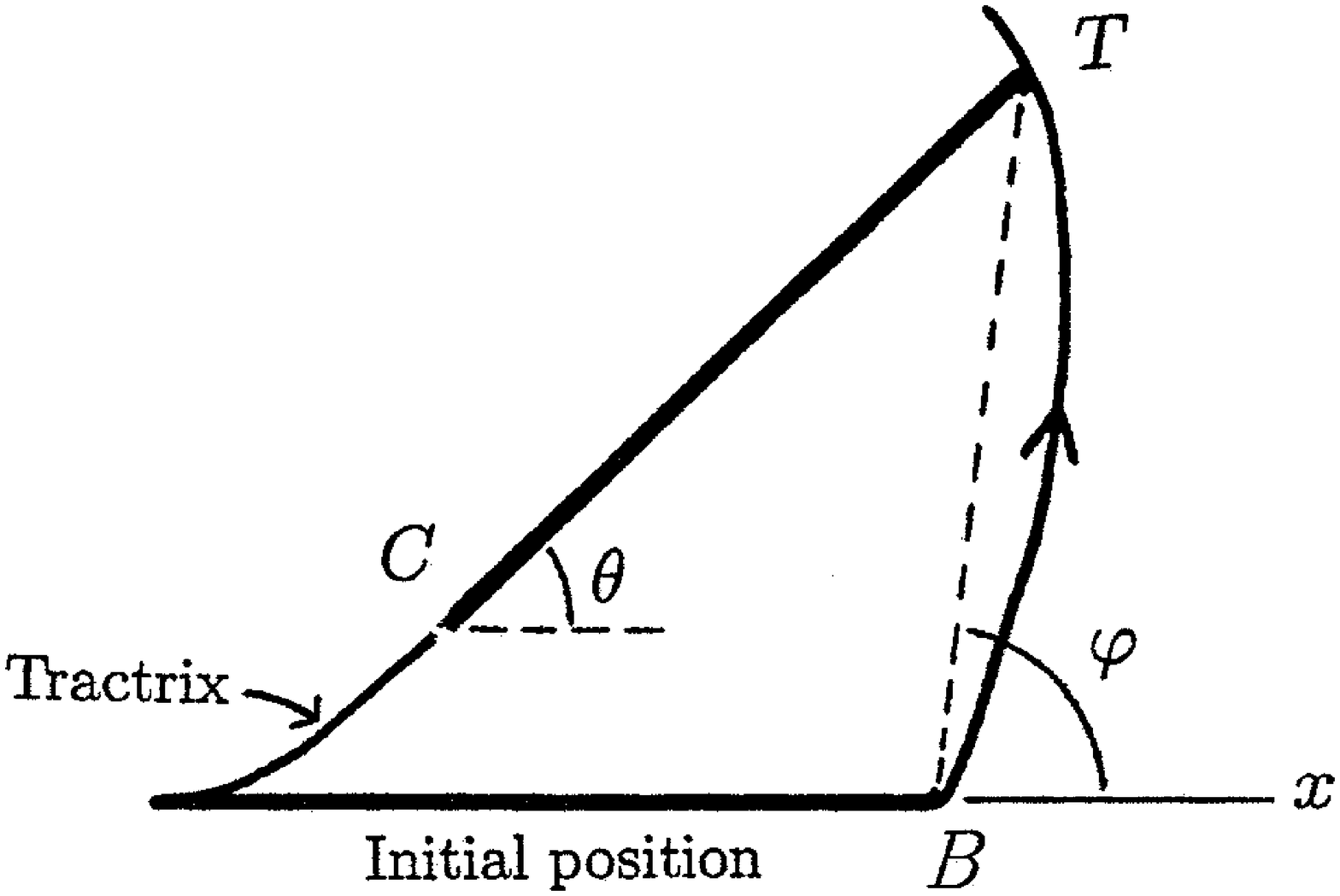}}
\centerline{Figure 8.}
\bigskip

%\vskip 0pt plus -1fill

%\newpage

Both authors use \thetag1 to write an infinite 
series expansion for the rate at which $\theta$ changes when the tracer point moves 
along an arbitrary curve (Figure~8). Hill, who keeps more terms, gets%\par}
$$
\l^2\der\theta\phi 
    \quad=\quad \(\frac{r^2}2 + \frac{r^4}{8\l^2} + \frac{r^6}{144\l^4} + \dots\)
    \quad+\quad \(\frac{r^3}{3\l} + \frac{r^5}{30\l^3} + \dots\)\cos(\phi - \theta).
$$
Here $(r,\phi)$ are polar coordinates of the tracer point $T$ about the base point $B$.

Their analysis depends on the function $r/\l$ being small on the path followed by the tracer 
point, which is to say, that $\O$
is small and that the base point $B$ is not far from $\O$, both relative to $\l$, the length of
the planimeter.
As the tracer point moves around $\bdyO$, the dominant terms not depending on $\theta$ integrate 
to
$$
\int_\bdyO \frac{r^2}2\,d\phi = A \qqand 
\int_\bdyO \frac{r^4}{8\l^2}\,d\phi = \frac1{2\l^2}\int_\O r^2\,dA = \frac{I_B}{2\l^2},
$$
where $A$ is the area of $\O$ and $I_B$ is the second moment (the moment of inertia) 
of $\O$ about $B$. 

Assuming the initial direction, $\theta=0$, of the planimeter is the positive $x$-axis, and that 
$|\theta|$ stays small (which follows from the assumption that $r/\l$ is small), Hill writes the 
remaining most dominant term as
$$
\frac1{3\l} \int_\bdyO{r^3}\cos(\phi - \theta)\,d\phi
   = \frac1{3\l} \int_\bdyO r^3 \(\cos\phi + \theta\sin\phi 
                                   - \frac12 \theta^2\cos\phi + \dots\)\,d\phi,
$$
and observes that 
$$
\frac1{3\l} \int_\bdyO r^3 \cos\phi \,d\phi 
    = \frac1\l \int_\O x\,dA = \frac{A\conj x}\l,
$$
where $\conj x$ is the $x$-coordinate of the centroid. Prytz and Hill thus have the approximation
$$
\l\s = \l^2\,\Delta\theta \approx A + \frac{A\conj x}\l + \frac{I_B}{2\l^2} 
     = A \(1 + \frac{\conj x}\l + \frac{R_B^2}{2\l^2}\),
$$
where $R_B^2$ is the mean-square distance of points in $\O$ from $B$. Although they don't say it
explicitly, it is easily seen that the error in this approximation is 
$\Cal O\big((d/\l)^3\big)$, where $d$ is the diameter of the set $\O \union \{B\}$ (or simply
the diameter of $\O$ if $B \in \conj\O$). Thus the error made by the approximation
$A \approx \l\s$ is $\Cal O(d/\l)$.  

Note that when measuring small regions,
$\Delta\theta$ is small but non-zero, that is, the chisel edge does not come to rest
in the place where it started. As we will see in the next section and at the end of
the last section, this can happen if the planimeter traces the boundary of a large 
region or a curve that bounds an oriented area of 0.

Prytz and Hill recommend measuring the region twice with the
same base point but with opposite initial directions, and averaging the results. Prytz's 
reasoning on this is not clear, but Hill notes that the $\conj x$ term will have opposite signs
in the two measurements, and so will drop out in the average, giving the approximation
$$
\l\conj\s \approx A + \frac{I_B}{2\l^2} = A \(1 + \frac{R_B^2}{2\l^2}\),
$$
where $\conj\s$ is the average of the displacements in the two measurements. From 
this, the error in the approximation $A \approx \l\conj\s$ is 
$\Cal O\big((d/\l)^2\big)$. Note that $\conj x$ also vanishes if the centroid is on 
the line through $B$ perpendicular to the initial position of the planimeter, agreeing
with the intuitive observation in Section~1 that a good starting position is one that 
is perpendicular to a line bisecting the region. 

The term $I_B = A R_B^2$ is minimized when $B$ is the centroid, which accounts for 
the recommendation of Prytz and Hill that the tracing start and end at the centroid.

From this analysis it appears that the holonomy $\Delta\theta$ of the Prytz planimeter
is some combination of {\sl all\/} of the moments of $\O$ about the base point $B$.
It is hoped that the modern approach taken in the next sections will lead to
a better understanding of the relationship between these.

%%% TXP %%% END Prytz123.tex %%%%%%%%%%%%%%%%%%%%%%%%%%%%%%%%%%%%%%%%%%%%%%

%%% TXP %%% BEGIN Prytz4.tex %%%%%%%%%%%%%%%%%%%%%%%%%%%%%%%%%%%%%%%%%%%%%%

\heading 
4. Motion of the Prytz Planimeter as Parallel Translation in a Fiber Bundle
\endheading

The configuration space for the Prytz planimeter is $E = \R^2 \cross S^1$, where the 
first factor is the location of the tracer point and the second factor is the angle 
giving the direction of the planimeter.  We will consider $E \To\pi \R^2$ as a trivial
circle bundle, where $\pi$ is projection onto the first factor.  If the tracer point
follows a curve $\g$ in $\R^2$, the resulting motion of the planimeter defines a
section of this bundle along $\g$.  We will see that this motion is described 
as parallel translation of an Ehresmann connection on $E \To\pi \R^2$. Thus, when 
the planimeter traverses the boundary of some region, its net rotation $\Delta\theta$ 
is an example of the holonomy of this connection. Our main source for 
connections and parallel translation in fiber bundles is \cite{KMS}, but also see 
\cite{Poo}. Many of the results in this and the next section illustrate the general 
theory in \cite{KMS}. 

Intuitively, the motion of the tracer point along a curve induces a one-parameter 
family of diffeomorphisms of the fiber $S^1$ (in which the fibers over different 
points are identified by projection $E\to S^1$ onto the second factor, \ie, by
the triviality of the bundle). It is clear that these 
diffeomorphisms are not simply rotations. For example, if the tracer point is moved 
along a straight line, the diffeomorphisms have two common antipodal fixed points, one
attracting and one repelling, the repelling fixed point being the fiber element in 
the direction of motion. As another example, consider the holonomy determined by
tracing the boundary of a region that is small relative to the length of the
planimeter. As seen in Sections~1 and~3, the resulting diffeomorphism
of the initial circle (which is the fiber over the base point where the tracing
begins and ends) is only approximately a rotation, since the net rotation of the
planimeter depends on its initial direction. Thus $E \To\pi \R^2$ should {\sl not\/} 
be viewed as a principal bundle in this context, since the group $G$ acting on the 
fiber is not the group of rotations (more precisely, the connection on 
$E \To\pi \R^2$, defined below, is not principal).  In the next section we will
determine the group $G$, and see that $E \To\pi \R^2$ is, in fact, an associated bundle
of a $G$-principal bundle, and that the connection on $E \To\pi \R^2$ is induced by
a connection on this principal bundle.

Take $(x,y,\theta)$ as coordinates on $E$, where $(x,y)$ are base coordinates and 
$\theta$ is the fiber coordinate.  The tracer point and chisel edge then have 
coordinates $q = (x,y)$ and $p = (x + \l\cos\theta, y + \l\sin\theta)$, respectively, 
and the forward-pointing normal (from Section~1) is 
$N = \sin\theta\pder{}x -\cos\theta\pder{}y$. 

%{\rightskip3truein %picture of coord sys
The expression $d\s = N{\cdot}dp$, defined in Section~1, is a 1-form on $E$. It is 
not exact, and henceforth will be denoted simply as $\s$.
The coordinate expression for $\s$ is %\par}
$$
\s = \sin\theta\,dx - \cos\theta\,dy - \l\,d\theta,
$$
and so the motion of the planimeter is governed by the non-holonomic condition
$$
\s = 0, \qquad\text{or}\qquad
\l\,d\theta = \sin\theta\,dx - \cos\theta\,dy.  \tag3
$$

For $e = (q,\theta) \in E$, let $H_e = \ker\s_e \subset T_eE$. The vector fields 
$$
X = \pder{}x + \frac1\l\sin\theta\pder{}\theta \qqand
Y = \pder{}y - \frac1\l\cos\theta\pder{}\theta         \tag4
$$
form a basis for $H_e$ at each $e \in E$, and 
$\pi_*|_{H_e} \: H_e \to T_{\pi(e)}\R^2$ is clearly an isomorphism, where $\pi_*$ is
the differential of $\pi$. It follows that $H = \{H_e\,\st\,e\in E\}$, viewed as a
sub-bundle of $TE \to E$, forms the distribution of horizontal subspaces of a 
connection on $E \To\pi \R^2$ (see \cite{KMS, \S9.3}). The connection form is the
1-form $\Phi = -\frac1\l \s \tensor \pder{}\theta$, which takes values in the bundle
$V$ of vertical vectors, where $V_e = \ker(\pi_*)_e \subset T_eE$. More precisely,
$\Phi_e \: T_eE \to V_e$, defined by 
$\Phi_e(W) = -\frac1\l \s(W)\left.\pder{}\theta\right|_e$, is projection onto $V_e$
with kernel $H_e$. If the tracer point of the planimeter follows a piecewise smooth
curve $\g$ in
$\R^2$, the resulting motion of the planimeter defines a curve $\gt$ in $E$ covering
$\g$ such that $\gt'(t) \in H_{\gt(t)}$ for all $t$, by virtue of \thetag3. Thus the 
motion of the planimeter is parallel translation for this connection.

The connection is also determined by its Christoffel form $\G$, which is a 1-form
on $\R^2$ with values in $\X$, the Lie algebra of smooth vector fields on $S^1$. 
(See \cite{KMS,\S9.7}. There needs to be just one Christoffel form since the bundle 
is trivial.) Suppose the tracer point of the planimeter is at $q\in\R^2$ moving with 
velocity $v = a\pder{}x + b\pder{}y$. The induced vector field 
$w = \G(v) = c(\theta)\pder{}\theta$ on $S^1$ is given by applying \thetag3:
$$
\l\, c(\theta) = \l\,d\theta(w\lowsub\theta) = (\sin\theta\,dx - \cos\theta\,dy)(v)
     = a\sin\theta - b\cos\theta.
$$
Consequently,
$$
\G(v)\lowsub\theta = w\lowsub\theta 
   = \frac1\l(a\sin\theta - b\cos\theta)\pder{}\theta. \tag5
$$
More directly, given $\theta \in S^1$, let 
$v\lowsub\theta = v + 0\pder{}\theta \in T_{(q,\theta)}E$.
Then $w = \G(v)$ is defined by
$\G(v)\lowsub\theta = -\Phi(v\lowsub\theta) = \frac1\l\s(v\lowsub\theta)\pder{}\theta$
(an identification is being made here between the fiber $E_q = \{q\}\cross S^1$ and the
model fiber $S^1$).  In either case we can write
$\G = \frac1\l(\sin\theta\,dx - \cos\theta\,dy)\tensor\pder{}\theta$.

In general, a connection on a fiber bundle need not be complete, that is, parallel
translation need not be defined for all elements of the bundle along all piecewise
smooth curves. However if the fiber is compact, as in the present case, then the
connection is necessarily complete, and so is properly called an Ehresmann connection.
This follows from a comment in \cite{KMS, \S9.9}, or by the following standard 
argument. Suppose $\g\:[0,1]\to\R^2$ is a piecewise smooth curve. Let $\epsilon>0$ and
extend $\g$ to have domain $I=(-2\epsilon,1+2\epsilon)$, and
consider the pull-back connection on the pull-back bundle $\g^*E \to I$. Let 
$f\:I\to\R$ be a function with support in $(-\epsilon,1+\epsilon)$ that is 
identically 1 on $[0,1]$. Let $T$ be the horizontal lift of the vector field
$f\der{}t$ to $\g^*E$. Then $T$ has compact support (by the compactness of the fiber),
and so is a complete vector field.  As the flow of $T$ over $[0,1]$ represents 
parallel translation along $\g$ from $\g(0)$ to $\g(1)$, the connection is complete. 
Note that this implies the intuitive observation made at the beginning of this section
that the motion of the tracer point along a curve 
induces a one-parameter family of diffeomorphisms of $S^1$, or equivalently, a
family of diffeomorphisms between the fibers of $E$ over points along the curve.

The vector fields $X$ and $Y$, given in \thetag4, are, respectively, the horizontal 
lifts of $\pder{}x$ and 
$\pder{}y$ for this connection. The curvature is given by the Lie bracket
$R(X,Y) = [X,Y] = \frac1{\l^2}\pder{}\theta$ (see \cite{KMS, \S9.4}). 
The vector fields $X$, $Y$, $\pder{}\theta$ form a frame on $E$. The dual frame is
easily seen to be $dx$, $dy$, $-\frac1\l\s$. In terms of this frame we have
$R = \frac1{\l^2}(dx\wedge dy)\tensor\pder{}\theta$. 
We see that if the tracer point moves 
counterclockwise around an infinitesimal rectangle with edges $\pder{}x = \pi_*(X)$ 
and $\pder{}y = \pi_*(Y)$, then the infinitesimal motion of the fiber is 
$\frac1{\l^2}\pder{}\theta$, which is also counterclockwise. It follows that the
holonomy group (which will be determined in the next section) must contain the rigid 
rotations, although from the results of
Sections~1 and ~3, it contains more than this.

Although the connection defined here on $E\to\R^2$ is not principal, it's interesting 
to observe that the curvature can be computed using a ``covariant exterior 
derivative'' as on a principal $S^1$-bundle \cite{KMS, \S11.5}. If $\a$ is a 1-form
on $E$, then $d_h\alpha$ is defined to be the horizontal component of $d\alpha$, that
is, $d_h\alpha(X,Y) = d\alpha(PX,PY)$ for $X,Y\in T_eE$, where $P\:T_eE\to T_eE$ is
projection onto $H_e$ with kernel $V_e$. Let $\o = -\frac1\l\s$ so that 
$\Phi = \o\tensor\pder{}\theta$. One computes that
$d\o = -\frac1\l\o\wedge(\cos\theta\,dx + \sin\theta\,dy) 
           - \frac1{\l^2}\,dx\wedge dy$.
The horizontal component of this is $\Omega = d_h\o = - \frac1{\l^2}\,dx\wedge dy$.
The curvature form $\Omega$ is related to the curvature tensor $R$ in the usual way,
namely, $R = -\Omega\tensor\pder{}\theta$.

%%% TXP %%% END Prytz4.tex %%%%%%%%%%%%%%%%%%%%%%%%%%%%%%%%%%%%%%%%%%%%%%%%

%%% TXP %%% BEGIN Prytz56.tex %%%%%%%%%%%%%%%%%%%%%%%%%%%%%%%%%%%%%%%%%%%%%

\heading 5. $E \To\pi \R^2$ as an Associated Bundle \endheading

We will now see that the bundle $E = \R^2 \cross S^1 \To\pi \R^2$ is an associated 
bundle of a principal bundle, and that the connection on $E$ is induced from a
connection on this principal bundle. As a comparison, and to further Bishop's analogy 
of the motion of the planimeter with that of a car, see \cite{Fe}. Our main references
for principal bundles are \cite{KMS, KN}.

As noted in the previous section, since the connection on $E\to\R^2$ is complete, 
motion of the tracer point along a piecewise $C^1$ curve in $\R^2$ induces a 
one-parameter family of diffeomorphisms.  The diffeomorphisms induced by all piecewise
$C^1$ curves lie in some smallest subgroup of $\Diff$, the group of all 
diffeomorphisms of $S^1$. Our first task is to determine this subgroup.

\proclaim{Theorem 1} The diffeomorphisms of $S^1$ induced by moving the tracer point
of the planimeter along arbitrary piecewise $C^1$ curves in $\R^2$ form a group,
namely, the group $\M$ of M\"obius transformations that preserve $S^1$ and its 
orientation.
\endproclaim

The proof consists of a number of steps. In the proof and the rest of the paper we
will write elements of $S^1$ as $e^{i\theta}$ and identify $\R^2$ with $\C$ when 
convenient.

First we show that the collection of diffeomorphisms is a group. The only question is
whether the collection is closed under composition. If $\g\:[0,1]\to\R^2$ is a 
piecewise $C^1$ curve, let $\psi_\g$ denote the diffeomorphism of $S^1$ induced by 
moving the tracer point along $\g$ from $\g(0)$ to $\g(1)$. Since the connection on 
$E\to\R^2$ is invariant under translations of $\R^2$, then $\psi_{\g+v} = \psi_\g$ for
every $v \in \R^2$. If $\g\lowsub1, \g\lowsub2\:[0,1]\to\R^2$ are two such curves, 
then $\psi_{\g\lowsub2}\comp\psi_{\g\lowsub1}$ is the diffeomorphism induced by moving
the tracer point along the curve $\g$ defined by $\g(t) = \g\lowsub1(2t)$ for 
$t\in[0,\frac12]$ and $\g(t) = \g\lowsub2(2t-1) + \g\lowsub1(1) - \g\lowsub2(0)$ for
$t\in[\frac12,1]$. Thus the collection of diffeomorphisms is closed under composition,
and hence is a group.

Before continuing with the proof of Theorem 1, we note some facts related to $\M$.
\roster
\item $\M$ is the same as the group of M\"obius transformations that preserve the 
      unit disk $D$, which is, of course, the group of orientation-preserving 
      isometries of the Poincar\'e disk model of the hyperbolic plane. These M\"obius
      transformations can be written as 
      $z \mapsto \frac{a z + b}{\conj b z + \conj a}$. This is a standard fact from
      elementary complex analysis, of course, but is done particularly nicely in
      \cite{Se}.
\item The ``mixed signature'' special unitary group
      $
      \SU = \left\{\smallpmatrix a & b \\ \conj b & \conj a \endsmallpmatrix
                      \: |a|^2 - |b|^2 = 1 \right\}
      $
      is a double cover of $\M$, the homomorphism being the usual map taking
      $A = \smallpmatrix a & b \\ \conj b & \conj a \endsmallpmatrix$ to the
      transformation $A{\cdot} z = \frac{a z + b}{\conj b z + \conj a}$. We will use 
      this to identify $\M$ with the projective group $P\!\SU = \SU/\{\pm I\}$, and 
      consequently to identify the Lie algebra of $\M$ with
      $\su = \left\{\smallpmatrix i\g & \b \\ \conj\b & -i\g \endsmallpmatrix
                                         \: \g\in\R,\, \b\in\C \right\}$.
      This is an easy exercise, and is done implicitly in \cite{Se}. It is completely
      analogous to the more common identification of $SL(2,\R)$ as a double cover of
      the M\"obius transformations preserving the upper half-plane (see, \eg, 
      \cite{Bo}). $\SU$ and $SL(2,\R)$ are conjugate subgroups of $SL(2,\C)$, related 
      by any matrix representing a M\"obius transformation mapping $D$ onto the upper 
      half-plane.
\item By restricting to $S^1$, $\M$ injects into $\Diff$. The differential of this is 
      a Lie algebra injection $\phi:\su\to\X$. More precisely, let $X\in\su$. To compute
      the vector field $\phi(X)$ at $z_0 = e^{i\theta_0} \in S^1$, let 
      $e^{i\theta(t)} = e^{tX}{\cdot}z_0$. Then
      $$
      \phi(X)_{z_0} = \theta'(0)\pderat{}\theta{z_0}
            = \({-}i\conj z_0 \derat{}t0 e^{tX}{\cdot}z_0\)\pderat{}\theta{z_0}.\tag6
      $$
\endroster

Next we determine how the Christoffel form $\G\:T\R^2\to\X$ and the map 
$\phi\:\su\to\X$ are related.

\proclaim{Lemma 2}  The Christoffel form $\G$ takes values in the image of $\phi$. In 
particular, if $v = a\pder{}x + b\pder{}y$, then $\G(v) = \phi(X)$, where
$X = -\frac1{2\l}\smallpmatrix 0 & a+ib \\ a-ib & 0 \endsmallpmatrix$.
\endproclaim

\demo{Proof} Suppose the tracer point is at $q\in\R^2$ moving with velocity
$v = a\pder{}x + b\pder{}y$. From \thetag5, the induced vector field on $S^1$ is
$$
\G(v)_{z_0} = \frac1\l(a\sin\theta_0 - b\cos\theta_0)\left.\pder{}\theta\right|_{z_0},
    \tag7
$$
where $z_0 = e^{i\theta_0}$. We need to identify $\G(v)$ as $\phi(X)$ for some 
$X \in \su$.

Let $X = \smallpmatrix 0 & \b \\ \conj\b & 0 \endsmallpmatrix$, and let
$\smallpmatrix a(t) & b(t) \\ \conj b(t) & \conj a(t) \endsmallpmatrix = e^{tX}$.
From \thetag6 we get $\phi(X)_{z_0} = 2\Im(\b \conj z_0)\pderat{}\theta{z_0}$.
If $\b = a + bi$, then 
$\phi(X)_{z_0} = 2(b\cos\theta_0 - a\sin\theta_0)\left.\pder{}\theta\right|_{z_0}$.
Comparing this with \thetag7, it follows that $\G(v) = \phi\big({-}\frac1{2\l}X\big)$.
\qed
\enddemo

Let $\g\:[0,1]\to\R^2$ be a $C^1$ curve. We seek a curve in $\M$ (actually its lift
in $\SU$) that generates the same one-parameter family of diffeomorphisms as $\g$. 
Suppose $A\:[0,1]\to\SU$ is a $C^1$ curve, and consider the resulting flow on $S^1$.
At $t=t_0$ the vector field for the flow is given (in a computation similar to that 
in \thetag6) by
$$
\({-}i\conj z_0 \derat{}t{t_0} A(t)A(t_0)\inv {\cdot}z_0\)\pderat{}\theta{z_0} 
      = \ \phi\big(A'(t_0)A(t_0)\inv\big)_{z_0}.
$$
For $A$ and $\g$ to generate the same vector field on $S^1$ at time $t$, we must have
$\G(\g'(t)) = \phi(A'(t)A(t)\inv)$.
Define the $\su$-valued 1-form $\o$ on $\R^2$ by $\o(v) = -\phi\inv(\G(v))$, 
that is, $\o(v)$ is the matrix in $\su$ that gives rise to the vector field $-\G(v)$ 
on $S^1$. From the lemma, we see that the formula for $\o$ is
$\o = \frac1{2\l}\pmatrix 0 & dz \\ d\conj z & 0 \endpmatrix$.
Then the desired curve in $SU(1,1)$ is the solution of the initial value problem
$$
A'(t) A(t)\inv = -\o(\g'(t)),\qquad A(0)=I. \tag8
$$

It follows that the group of diffeomorphisms of $S^1$ generated by the motion of 
the tracer point is a subgroup of $\M$.  By the lemma, the values of $\o$ are of the 
form $\smallpmatrix 0 & \b \\ \conj\b & 0 \endsmallpmatrix$. These generate the 
entire Lie algebra $\su$, and so it follows that the group is $\M$. This completes the
proof of Theorem~1.

\medskip

From the proof of the theorem we see that the action of the one-parameter subgroup 
$e^{-t\o(v)}$ on $S^1$ is the same as that of moving the tracer point in a straight 
line in $\R^2$ with constant velocity $v$. If the tracer point starts at $q \in \R^2$ 
with the planimeter in the direction $e^{i\theta_0}$, then the resulting curve in 
$E = \R^2\cross S^1$ is $(q + tv, e^{-t\o(v)}{\cdot} e^{i\theta_0})$. Identifying
$\C$ and $\R^2$, the path of the chisel, which is a standard tractrix, is 
$\tau(t) = (q + tv) + \l e^{-t\o(v)}{\cdot} e^{i\theta_0}$. Letting $q=0$ and 
$v = \pder{}x$ (identified with $1 = 1 + 0i$ in $q+tv$), we have
$e^{-t\o(v)} = \smallpmatrix \cosh(t/2\l) & -\sinh(t/2\l) \\
                            -\sinh(t/2\l) &  \cosh(t/2\l)\endsmallpmatrix$ and
$$
\align
\tau(t) &= t + \frac{\cosh(t/2\l)e^{i\theta_0} - \sinh(t/2\l)}
                    {-\sinh(t/2\l)e^{i\theta_0} + \cosh(t/2\l)} \\
        &= t + \l\frac{\cosh(t/\l)\cos\theta_0  - \sinh(t/\l)}
                      {\cosh(t/\l) - \sinh(t/\l)\cos\theta_0}
             + i\l\frac{\sin\theta_0}{\cosh(t/\l) - \sinh(t/\l)\cos\theta_0}.
\endalign
$$
When $\theta_0 = \pi/2$, this simplifies to 
$\tau(t) = t - \l\tanh(t/\l) + i\l\sech(t/\l)$, which is \thetag2 in complex form.

Now consider $S^1$ as the points at infinity of the Poincar\'e disk model $D$ of the 
hyperbolic plane. Acting on $D$, $e^{-t\o(v)}$ is the one-parameter group of 
hyperbolic translations that moves the origin along a geodesic with initial velocity 
$-\frac v{2\l}$ (note that this is in the opposite direction of the motion of the 
tracer point). This has the following nice interpretation. Suppose that a hyperbolic 
stargazer goes walking with constant velocity $v$ in $D$.  It is natural for the 
stargazer to think of herself as always being at the center of the celestial circle, 
with the hyperbolic plane passing beneath her feet with constant velocity $-v$. If she
fixes her gaze on a particular star, then the retrograde motion of that star is the 
same as the motion of a (Euclidean!)\ Prytz planimeter of length $\l = 1/2$.

This discussion shows how a geodesic in $\R^2$ can be developed into a geodesic in 
$D$. Generalizing this, one can develop any polygonal path in $\R^2$ to a polygonal 
path in $D$ (and by polygonal approximations, any piecewise $C^1$ path). 
One needs to be careful, however. A sequence $e^{X_1}$,\dots, $e^{X_n}$ of 
translations that moves the origin around a closed path in $D$ determines a rotation 
(relative to the origin)
of the points at infinity, whereas the tracer point of the planimeter following a 
closed loop typically does not result in a pure rotation of the initial directions of 
the planimeter. Evidently, the motion of the tracer point that induces the same 
sequence $e^{X_1}$,\dots, $e^{X_n}$ is {\sl not\/} generally a closed path in $\R^2$.
Similarly, the motion of the tracer point in $\R^2$ around a closed path does not
generally develop to a closed path in $D$.

As a consequence of Theorem~1, we view $E\to\R^2$ as an $\M$-bundle or as an
$\SU$-bundle \cite{KMS, \S10.1}. The form $\o$ in the proof is, of course, a ``local 
frame representation'' of a connection on some principal bundle \cite{KN, pg.~66}
(a ``physicist's connection'' in \cite{KMS, \S11.4}), and \thetag8 is the 
corresponding parallel translation equation. Writing down the principal bundle and 
realizing $E\to\R^2$ as an associated bundle is now straight forward.

For the remainder of this section, let $G$ be $\M$ or $\SU$. Consider the principal 
$G$-bundle $P = \R^2 \cross G \to \R^2$. As 
$G$ acts on $S^1$, we have the standard construction of the associated $S^1$ 
fiber bundle $\tilde E \to \R^2$ (see \cite{KN, pg.~54; KMS, \S10.7}). The space 
$\tilde E$ is given by $(P \cross S^1)/{\sim}$, in which 
$(q,A,e^{i\theta}) \sim (q,I,A{\cdot}e^{i\theta})$. The map $E \to \tilde E$ that 
takes $(q,e^{i\theta})$ to the equivalence class of $(q,I,e^{i\theta})$ is clearly a
fiber bundle equivalence. 

Each element $(q,A)$ in the fiber $\{q\}\cross G$ of $P \to \R^2$ represents a 
``frame'' for the fiber $\{q\}\cross S^1$ of $E \to \R^2$, that is, a diffeomorphism 
from the model fiber $S^1$ to $\{q\}\cross S^1$ given by 
$e^{i\theta} \mapsto (q,A{\cdot}e^{i\theta})$. As both bundles are trivial, it is 
easiest to make computations relative to the ``standard frame,'' that is, the 
``identity'' diffeomorphism given by projection $\{q\}\cross S^1 \to S^1$, represented 
by $(q,I) \in P$. The standard frame is more than a computational convenience, 
however, since it represents the Euclidean geometry of the plane. The identification
it makes of the fibers of $E\to\R^2$ is by Euclidean translation. Writing
the connections on $E\to\R^2$ and $P\to\R^2$ using the ``local descriptions'' of
$\G$ and $\o$ amounts to describing how the motion of the planimeter differs from
Euclidean parallel translation of vectors.

The connection form $\ot\:TP\to\su$ along the identity section of $P \to \R^2$ is 
given as follows. If $v{\dirsum}X \in T_{(q,I)}P = T_{q}\R^2 \dirsum \su$, then 
$\ot(v{\dirsum}X) = \o(v) + X$. Note that $\ot(v{\dirsum}X) = 0$ implies $X = -\o(v)$,
that is, $v{\dirsum}X$ is horizontal when $X$ induces the same vector field on $S^1$ 
as the planimeter when the tracer point undergoes the infinitesimal displacement $v$. 
The connection form $\ot$ is extended off of the identity section by the usual 
equivariance requirement \cite{KN, pg.~64; KMS, \S11.1}: if 
$v{\dirsum}XA \in T_{(q,A)}P$, then 
$\ot(v{\dirsum}XA) = \Ad(A\inv)(\ot(v{\dirsum}X)) = \Ad(A\inv)(\o(v) + X)$.
It is clear from the proof of Theorem~1 that the induced connection on 
$E \to \R^2$ viewed as an associated bundle is the same as the connection 
$-\frac1\l\sigma\tensor\pder{}\theta$ in the previous section, as they induce the 
same parallel translation. This also follows from Lemma~2 and the theorem in 
\cite{KMS, \S11.9}.

Finally, we compute the curvature and holonomy. Relative to the standard frame, the 
curvature form is
$$
\Omega = d\o + \o \wedge \o
  = \frac1{4\l^2}\pmatrix dz \wedge d\conj z & 0\\ 0 & d\conj z \wedge dz \endpmatrix
  = -\frac1{\l^2}\pmatrix i/2 & 0\\ 0 & -i/2 \endpmatrix dx \wedge dy.
$$
The infinitesimal rotations in $\su$ are of the form 
$\smallpmatrix i\g & 0 \\ 0 & -i\g \endsmallpmatrix$. Thus the curvature is 
purely rotational, at least relative to the standard frame. If 
$Z = \smallpmatrix i\g & 0 \\ 0 & -i\g \endsmallpmatrix$, then
$e^Z = \smallpmatrix e^{i\g} & 0 \\ 0 & e^{-i\g} \endsmallpmatrix$,
which acts as a rotation through angle $2\g$. Remembering that $-\Omega$
should be a measure of infinitesimal holonomy, let
$Z = -\Omega\(\pder{}x,\pder{}y\) 
   = \frac1{\l^2}\smallpmatrix i/2 & 0\\ 0 & -i/2 \endsmallpmatrix$.
Then $e^{tZ}$ induces the infinitesimal rotation $\phi(Z) = \frac1{l^2}\pder{}\theta$, 
agreeing with the bracket curvature
computation $R(X,Y)=[X,Y]$ from the previous section, where $X$ and $Y$ are
the horizontal lifts in $E$ of $\pder{}x$ and $\pder{}y$ given in \thetag4.

For $p\in\R^2$ let $\Cal H_p \subset G$ be the holonomy group based at $p$ for the 
connection on $P\to\R^2$, and let $\frak h_p$ be its Lie algebra. The curvature 
computation shows that $\frak h_p$ contains the infinitesimal rotations of 
$E_p = \{p\}\cross S^1$, and so $\Cal H_p$ contains the rotations of $E_p$ (or their 
representations in $\SU$), even though we have not yet seen a closed curve that 
induces a rotation as its holonomy. By the Ambrose-Singer Theorem \cite{KN, pg.~89}, 
if $Z = \O\(\pderat{}xq,\pderat{}{\botsmash y}q\)$ is an infinitesimal rotation at 
some other 
point $q$, then $\psi_\g\inv Z \psi_\g$ is in $\frak h_p$, where $\psi_\g$ is parallel
translation for the connection in $P\to\R^2$ along some curve $\g$ from $p$ to $q$.
In particular, we can parallel translate along the segment joining $p$ and $q$.
Let $v = q-p$. Then the proof of Theorem~1 (particularly Lemma~2), shows that
parallel translation from $p$ to $q$ is represented by $e^X$ in $G$, where 
$X = -\o(v)$. Thus $e^{-X}Z e^X$ is in $\frak h_p$; it is the infinitesimal holonomy 
resulting from moving the tracer point along the segment from $p$ to $q$, around a 
small loop, and back along the segment to $p$. It follows that 
$\derat{}t0(e^{-tX}Z e^{tX}) = -[X,Z]$ is in $\frak h_p$. As we have seen, $X$ has the
form $\smallpmatrix 0 & \b \\ \conj\b & 0 \endsmallpmatrix$. But
$-[X,Z] = 2\smallpmatrix 0 & i\b \\ -i\conj\b & 0 \endsmallpmatrix$ has the same form,
and these generate the Lie algebra $\su$, as noted earlier. It
follows that $\Cal H_p$ is the entire group $G$. When $\Cal H_p$ acts on $E_p$, it 
follows that the holonomy group at $p$ for the bundle $E\to\R^2$ is $\M$. 

The holonomy for $E\to\R^2$ could also be determined by computing the Lie algebra
generated by the vector fields $\G(v)$, and appealing to the generalized 
Ambrose-Singer Theorem due to Michnor \cite{KMS, \S9.11}. An easy computation shows
that this Lie algebra is the image of the map $\phi:\su\to\X$ used in Lemma~2, which 
is isomorphic to $\su$.

This is all summarized in the following theorem.

\proclaim{Theorem 3} Let $G$ be $\M$ or $\SU$.
\roster
\item The bundle $E = \R^2\cross S^1 \to \R^2$ is (isomorphic to) an associated bundle
      of the principal bundle $P = \R^2\cross G\to\R^2$.
\item The motion of the planimeter is given by parallel translation for a connection 
      on $E \to \R^2$, which is induced by a principal connection on $P \to \R^2$.  
      With respect to the standard frame, the connection and curvature forms for the
      principal connection are 
      \,$\o = \frac1{2\l}\smallpmatrix 0 & dz \\ d\conj z & 0 \endsmallpmatrix$ and
      \,$\Omega = -\frac1{\l^2}\pmatrix i/2 & 0\\ 0 & -i/2 \endpmatrix dx \wedge dy$.
\item The holonomy group for the connection on $P \to \R^2$ is $G$.
\item The holonomy group for the connection on $E \to \R^2$ is $\M$.
\endroster
\endproclaim

\subheading{Questions and Speculations}
The motion of the tracer point around a closed curve, starting and stopping at a base
point, determines an element of the holonomy group. If a different base point 
on the curve is used, the two holonomy elements need not be the same, but are 
conjugate to each other. (The holonomy groups at the two different base points are 
identified by the Euclidean translation that identifies the corresponding fibers of 
$E\to\R^2$.) Thus a region $\O$ with a simple, closed, piecewise $C^1$ boundary 
determines a conjugacy class in $G$.  What information about $\O$ can be determined 
from its conjugacy class?

Suppose $\O$ is sufficiently small relative to the length $\l$ of the planimeter so
that the holonomy $H$ determined by tracing $\bdyO$ starting at base point $B$ acts
on $S^1$ without fixed points. When $H$ acts on the Poincar\'e disk $D$, it is a
hyperbolic rotation about some point $z(B) \in D$ (see the next section). Every 
hyperbolic rotation is conjugate to a rotation about the origin. In particular, let 
$\phi\lowsub B(\z) = \frac{\z - z(B)}{1 - \conj z(B)\z}$. Then $\phi\lowsub B$ is the 
hyperbolic translation taking $z(B)$ to the origin, and 
$\phi\lowsub B H\phi_B\inv$ is a rotation about the origin. The transformation 
$\phi\lowsub B$ is represented by 
$\smallpmatrix a & b \\ \conj b & a \endsmallpmatrix$ in $\SU$, where
$a = 1/\sqrt{1 - |z(B)|^2}$ and $b = -z(B)/\sqrt{1 - |z(B)|^2}$. From our results
above, this matrix is $e^{-\o(v)}$ for some $v\in\R^2$. Letting $f(B)$ be the point
$B+v$, it follows that $f(B)$ is the unique point in $\R^2$ with the following 
property. If the planimeter starts tracing at $f(B)$, moves along the segment joining
$f(B)$ to $B$, goes around $\bdyO$, and then back along the segment to $f(B)$, then 
the resulting holonomy is purely rotational. In this
way the curve $\bdyO$ and the number $\l>0$ determine curves $z\:\bdyO\to D$ and 
$f\:\bdyO\to\R^2$. Are these curves related to $\bdyO$ in some simpler way?

Given an element of the  holonomy group, it would be nice to find a loop, perhaps of 
shortest length, that generates it.  In particular, what curves generate pure 
rotations in $G$? Some are given in the previous paragraph, but can they be 
characterized more simply? What non-trivial curves generate the identity? When 
$G = \SU$, what curves generate holonomy $-I$?

For the Prytz planimeter on the plane, it's easiest to refer everything to the 
standard frame, since the bundles are all trivial. The full machinery of the principal
bundle viewpoint may be necessary to study how the planimeter works on the sphere.  
In this case the configuration space $E \to S^2$ is the unit circle bundle over the 
sphere, which would appear to be an associated bundle of an appropriate principal 
$G$-bundle over $S^2$, as in the $\R^2$ case. The Hopf fibration factors as 
$S^3 \to E \to S^2$ (the first map is a double cover), and so it seems that the 
connection on $E \to S^2$ given by the motion of the planimeter should lift to a 
connection on the Hopf fibration, which would have $\SU$ as its group.

\heading 6. A Special Case of Menzin's Conjecture \endheading

Recall Menzin's conjecture from Section~2: If a planimeter of length $\l$ traces the
boundary of a region $\O$ with area $A > \pi\l^2$, then the chisel edge asymptotically
approaches a particular tractrix that is closed. When the planimeter follows this 
closed tractrix, it makes a complete rotation. In this section we prove this in the 
case when $\bdy\O$ is a parallelogram.

The conjecture can be rephrased in terms of the Poincar\'e return map on the initial 
circle and the winding number of a section of $E \to \R^2$ over $\bdy\O$. Let $B$ be 
a base point on $\bdy\O$, and let $H\:S^1\to S^1$ be the diffeomorphism in $\M$ 
defined by tracing $\bdy\O$ in the positive direction with the planimeter starting 
and stopping at $B$, that is, $H$ is the holonomy defined by the curve 
$\bdy\O$ and the base point $B$. If the planimeter starts with a fixed point of $H$ 
as its initial direction, then the motion of the planimeter determines a section of 
$E \to \R^2$ over $\bdy\O$. Menzin's conjecture then becomes the following.

\proclaim{Menzin's Conjecture} If $A > \pi\l^2$, then $H$ has a unique attracting fixed point.
If the planimeter starts with this fixed point as its initial position, then the induced section 
of $E \to \R^2$ over $\bdy\O$ has winding number 1.
\endproclaim

Note that for topological reasons, one generically expects to get a repelling fixed 
point as well. Menzin's conjecture also predicts this---the repelling fixed point is 
the attracting fixed point when the tracing direction is reversed.

In contrast, if the planimeter traces a region that is small relative to its length, 
the results of Prytz and Hill imply that $H$ has no fixed points, as observed earlier.

We recall a standard fact about elements of $\M$ (actually their representatives in 
$\SU$), their fixed points on $S^1$, and how they act as isometries on the Poincar\'e 
disk $D$ (see, \eg, \cite{Se}): 
\roster
\item $A\in\SU$ has two fixed points on $S^1$ if and only if $\left|\tr A\right| > 2$.
      In this case $A$ is a hyperbolic translation in $D$ along the geodesic joining 
      the fixed points on $S^1$. One of the fixed points is attracting and the other 
      is repelling.  
\item $A$ has one fixed point on $S^1$ if and only if $\left|\tr A\right| = 2$ and 
      $A \ne \pm I$.
      In this case $A$ is a horocyclic rotation on $D$ (a motion of $D$ that preserves
      every horocycle through the fixed point). The fixed point is semi-attracting, 
      that is, it is attracting from one side, but repelling from the other. 
\item $A$ has no fixed points on $S^1$ if and only if $\left|\tr A\right| < 2$.
      In this case $A$ has a fixed point in $D$ and is a hyperbolic rotation about 
      this point.
\endroster
Evidently any proof of Menzin's conjecture 
will involve showing that $A > \pi\l^2$ implies $\left|\tr H\right| > 2$.

The results of Prytz and Hill combined with Menzin's conjecture yield qualitatively
different holonomies when tracing the boundaries of small and large regions. The
holonomy of the null curve (the curve that stays at the base point) is, of course,
$I$, which has trace 2, and so is on the boundary between the open sets $U_0$ and 
$U_2$ in $\SU$ of transformations with no fixed points and two fixed points. 
Traversing a small loop will yield a holonomy $H$ close to $I$.  If the loop is the
boundary of a region, then $H \in U_0$.  It is possible to traverse a small loop so
that $H \in U_2$, but the loop will necessarily have to enclose some area in the
positive sense and some in the negative sense, as in a figure eight. An example of
this will follow the proof of the special case.

\demo{Proof of Menzin's conjecture when $\bdy\O$ is a parallelogram} Let $v$ and $w$ 
form a positively oriented basis
of $\R^2$, and let $\bdy\O$ be the parallelogram with vertices 0, $v$, $v+w$, $w$. We will
compute the holonomy determined by tracing $\bdy\O$ in the positive direction with 0 as the
base point. Let
$$
\gather
X = \smallpmatrix 0 & \b \\ \conj\b & 0 \endsmallpmatrix = -\o(v) 
  = -\tfrac1{2\l}\!\smallpmatrix 0 & dz(v) \\ d\conj z(v) & 0 \endsmallpmatrix \tag9\\
\qqand
Y = \smallpmatrix 0\vphantom\b & \d \\ \conj\d & 0\vphantom\b \endsmallpmatrix = -\o(w) 
  = -\tfrac1{2\l}\!\smallpmatrix 0 & dz(w) \\ d\conj z(w) & 0 \endsmallpmatrix\!,
\endgather
$$
and note that 
$$
A = \Im(\conj v w) = 4\l^2\Im(\conj\b \d) = 4\l^2|\b\d|\sin\theta, \tag9
$$
where $\theta$ is the angle between $v$ and $w$.

We have 
$e^X = \smallpmatrix a & b \\ \conj b & a \endsmallpmatrix 
     = \smallpmatrix a & \a\b \\ \a\conj\b & a \endsmallpmatrix$ 
and
$e^Y = \smallpmatrix c & d \\ \conj d & c \endsmallpmatrix 
     = \smallpmatrix c & \g\d \\ \g\conj\d & c \endsmallpmatrix$, 
where $a = \cosh|\b|$, $\a = \frac{\sinh|\b|}{|\b|}$, $b = \a\b$,
$c = \cosh|\d|$, $d = \g\d$, and $\g = \frac{\sinh|\d|}{|\d|}$.

The holonomy is then $H = e^{-Y}e^{-X}e^Y e^X$, which can be written as
$$
H \ = \ I + \tilde H 
  \ = \ I + 2\Im(\conj bd)\pmatrix i & 0 \\ 0 & -i \endpmatrix e^Y e^X.
$$
Although it's not essential to the proof, it's nice to observe that
$$
\tilde H e^{-X}e^{-Y} =
2\Im(\conj bd)\pmatrix i & 0 \\ 0 & -i \endpmatrix 
   = \a\g \frac{A}{\l^2}\pmatrix i/2 & 0 \\ 0 & -i/2 \endpmatrix = \a\g [Y,X],
$$
which should be compared with the expressions for infinitesimal holonomy in the 
previous two sections.

We need to show that $\left|\tr H\right| > 2$. One easily computes that 
$\tr H = 2 - 4\Im^2(\conj bd)$. Thus $\tr H$ cannot be bigger than 2, and
$\tr H < -2$ if and only if $\Im(\conj bd) > 1$. It follows that $H$ has an attracting
fixed point if and only if
$$
\Im(\conj bd) = \a\g\Im(\conj\b \d) = \a\g|\b\d|\sin\theta = \a\g\frac A{4\l^2} > 1, 
\qquad\text{\ie,}\qquad A > \frac4{\a\g}\l^2. \tag{10}
$$
As $\a$ and $\g$ are larger than 1, this is close to the hypothesis $A > \pi\l^2$! 
Using the expressions for $\a$ and $\g$, we have that $H$ has an attracting fixed 
point if and only if $\,\Im(\conj bd) = \sinh|\b|\,\sinh|\d|\,\sin\theta > 1$.

Using \thetag9, the implication
$$
A > \pi\l^2 \qquad \implies \qquad \text{$H$ has an attracting fixed point}
$$ 
becomes
$$
|\b\d|\sin\theta > \frac\pi4 \qquad \implies \qquad 
            \sinh|\b|\,\sinh|\d|\,\sin\theta > 1.
$$
One easily finds that the minimum of \,$\sinh x\,\sinh y\,\sin\theta$ \, subject to 
the constraint
$xy\sin\theta \ge \pi/4$ is $(\cosh\sqrt\pi - 1)/2$, which is approximately 1.014, and so this
implication holds.

To see that the planimeter makes a full rotation as the chisel edge follows one of
the periodic trajectories, we compute the fixed points of $H$ and observe that the 
planimeter always rotates counterclockwise. If this is to 
happen, then the planimeter should make half a rotation as the tracer point moves from
one vertex of the parallelogram to the opposite vertex. Consequently we look for
solutions of $e^Y e^X{\cdot}z = -z$, that is,
$$
\frac{(ac + \conj bd)z + (ad + bc)}{(a\conj d + \conj bc)z + (ac + b\conj d)} = -z
$$
(this equation is considerably easier than $H{\cdot}z = z$).
Using $a^2 - |b|^2 = c^2 - |d|^2 = 1$, the discriminant of this quadratic simplifies
to $4(1 - \Im^2(\conj bd))$, which is negative by \thetag{10}. The solutions can then
be written as
$$
z_\pm = -\frac{ad + bc}{|ad + bc|}\,
            \frac{ac + \Re(\conj bd) 
               \pm i\sqrt{\botsmash{\Im^2(\conj bd) - 1}}}{|ad + bc|}
$$
(this is to be taken as $z_+$ and $z_-$, using the corresponding sign on the radical).
It is easily shown that these are indeed the fixed points of $H$.
Both fractions in this expression have unit modulus. Since
$a^2 c^2 = (|b|^2 + 1)(|d|^2 + 1) > |\conj bd|^2 \ge \Re^2(\conj bd)$, we have
$ac + \Re(\conj bd) > 0$, and so the two values
of the second fraction form a conjugate pair with positive real part. Thus the fixed
points of $H$ consist of this conjugate pair rotated through the angle
$\arg(-(ad + bc))$. Thought of as vectors in $\R^2$, $b$ and $d$ point in the 
directions {\sl opposite\/} those of the vectors $v$ and $w$ that give the sides of
the parallelogram (see \thetag9).  Thus it is plausible that the conjugate pair of the second 
fraction above are rotated into the interior of the angle formed by $v$ and $w$
(see Figure~9). This is in fact the case, as will now be shown.

%\vskip4.5truecm
\bigskip
\centerline{\epsfxsize.45\hsize\epsfbox{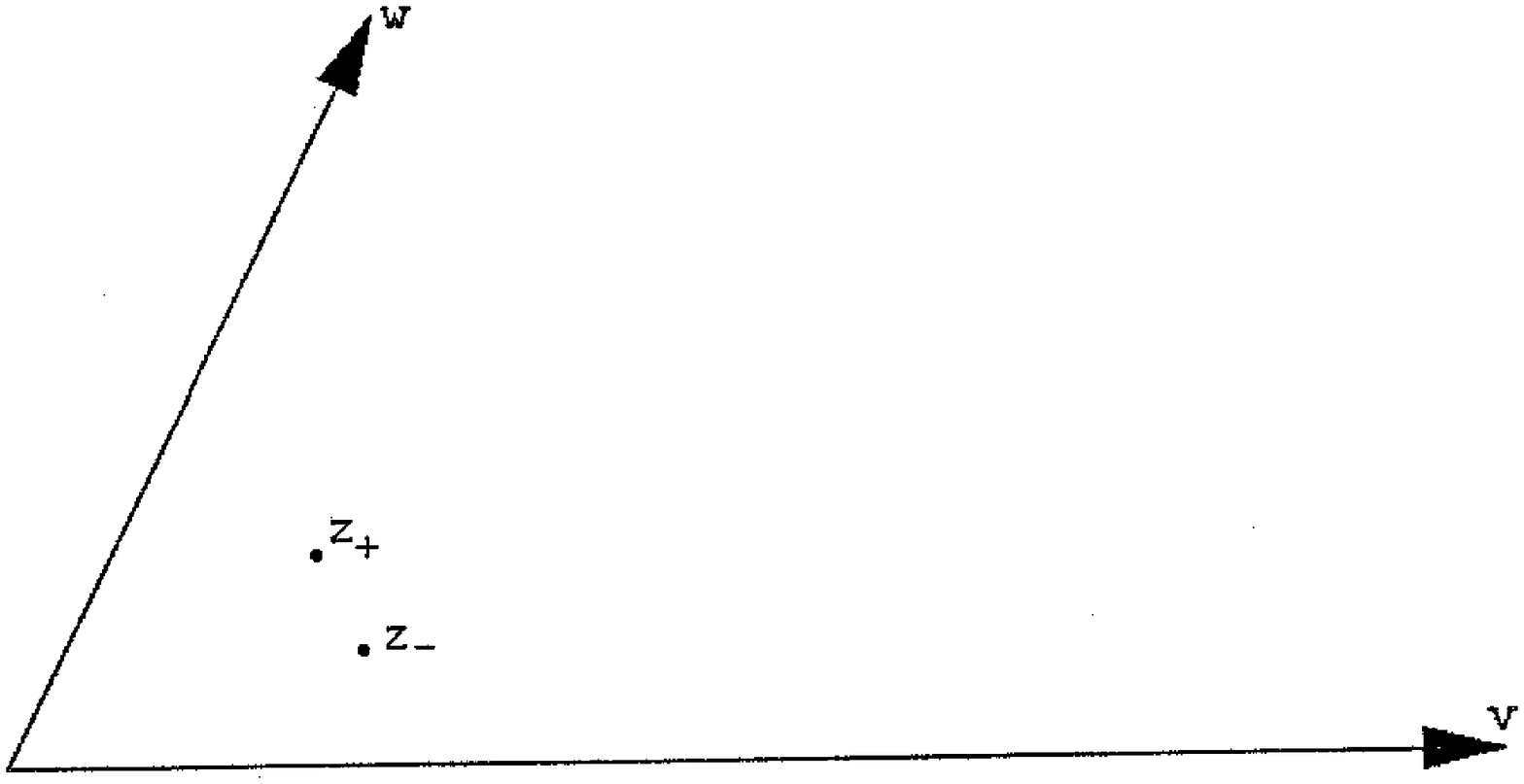}\qquad\epsfxsize.45\hsize\epsfbox{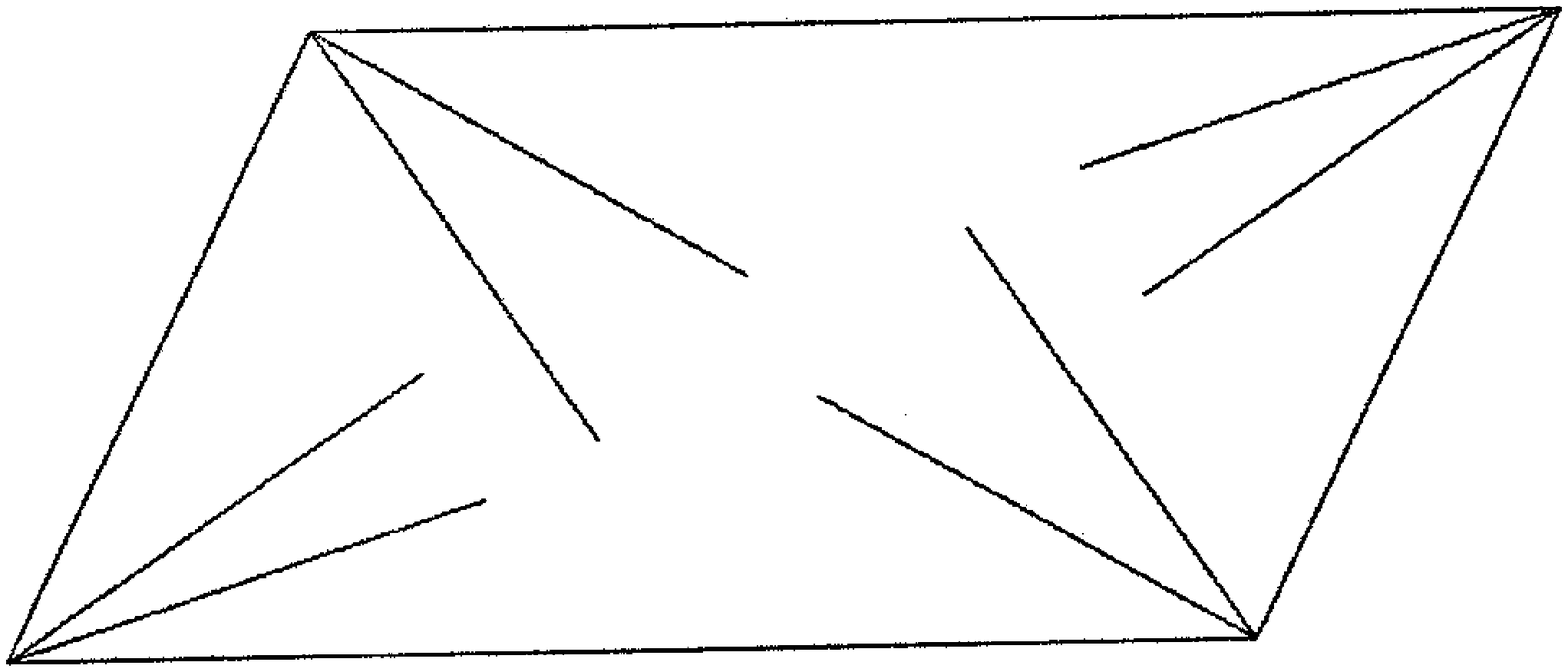}}
\centerline{Figure 9. \hfil\hfil  Figure 10.}
\bigskip

To be more precise, we need to show that $\arg(z_\pm/v) = \arg(-z_\pm/b)$ and
$\arg(w/z_\pm) = \arg(-d/z_\pm)$ are positive, thinking of $v$ and $w$ as complex
numbers. The condition \linebreak 
$\arg(-z_\pm/b)>0$ is equivalent to $\Im(-z_\pm\conj b)>0$, 
which we now show. Proving $\arg(-d/z_\pm)>0$ is similar.  We have
$$
|ad + bc|^2 \Im(-z_\pm\conj b) 
  = a\Im(\conj bd)\(ac + \Re(\conj bd)\)
     \pm \(a\Re(\conj bd) + |b|^2 c\)\sqrt{\Im^2(\conj bd) - 1}.
$$
For this to be positive, we need
$$
a\Im(\conj bd)\(ac + \Re(\conj bd)\) 
   > \left|a\Re(\conj bd) + |b|^2 c\right| \sqrt{\Im^2(\conj bd) - 1}.
$$
Since $\Im(\conj bd) > 0$ by \thetag{10} and $ac + \Re(\conj bd) > 0$ (as observed above),
this is equivalent to
$$
a^2\Im^2(\conj bd) \(ac + \Re(\conj bd)\)^2
   > \( a^2\Re^2(\conj bd) + 2a|b|^2 c \Re(\conj bd) + |b|^4 c^2 \) 
       \(\Im^2(\conj bd) - 1\).
$$
On the right hand side replace $a^2$ with $|b|^2 + 1$, and $|b|^4$ with 
$|b|^2 (a^2 - 1)$. Collecting the terms containing $|b|^2$ other than $|b|^2 c^2$, we 
get the equivalent inequality
$$
a^2\Im^2(\conj bd) \(ac + \Re(\conj bd)\)^2
   > \( |b|^2 \(ac + \Re(\conj bd)\)^2 + \Re^2(\conj bd) - |b|^2 c^2 \) 
       \(\Im^2(\conj bd) - 1\).
$$
Collecting the terms with $\(ac + \Re(\conj bd)\)^2$ and using $a^2 - |b|^2 = 1$, this
becomes
$$
\(|b|^2 + \Im^2(\conj bd)\) \(ac + \Re(\conj bd)\)^2
   > \( \Re^2(\conj bd) - |b|^2 c^2 \) \(\Im^2(\conj bd) - 1\).
$$
Using $c^2 = |d|^2 + 1$, we have 
$\Re^2(\conj bd) - |b|^2 c^2 = -\(|b|^2 + \Im^2(\conj bd)\)$. Thus all of these
inequalities are equivalent to
$$
\(|b|^2 + \Im^2(\conj bd)\) \(\(ac + \Re(\conj bd)\)^2 + \(\Im^2(\conj bd) - 1\)\) > 0.
$$
This inequality holds since $\Im^2(\conj bd) > 1$ by \thetag{10}, and so we have
$\arg(z_\pm/v) > 0$, as desired.

This argument shows that if the planimeter is following either of the periodic 
trajectories,
when the tracer point is at the initial vertex of the parallelogram the chisel edge
is in the interior of the angle formed by the adjacent edges.  It holds for the other
vertices as well. Figure~10 illustrates this situation.  As the tracer point moves 
along the initial edge of the parallelogram from 0 to $v$, the chisel edge follows a 
standard tractrix. It stays to the left of the direction of motion, that is, on the 
same side of the edge as the parallelogram, and thus only rotates counterclockwise. 
When the tracer point reaches $v$, the chisel edge is in the interior of the angle at
$v$, and so is to the left of the new direction of motion when the tracer point starts
moving from $v$ to $v+w$ (see Figure~11). This continues around the parallelogram, and
so the planimeter has made a full rotation when the tracer point returns to the 
origin. \qed
\enddemo

%\newpage
%\ \vfill
%\ \vskip4.5truecm
\centerline{\epsfxsize.7\hsize\epsfbox{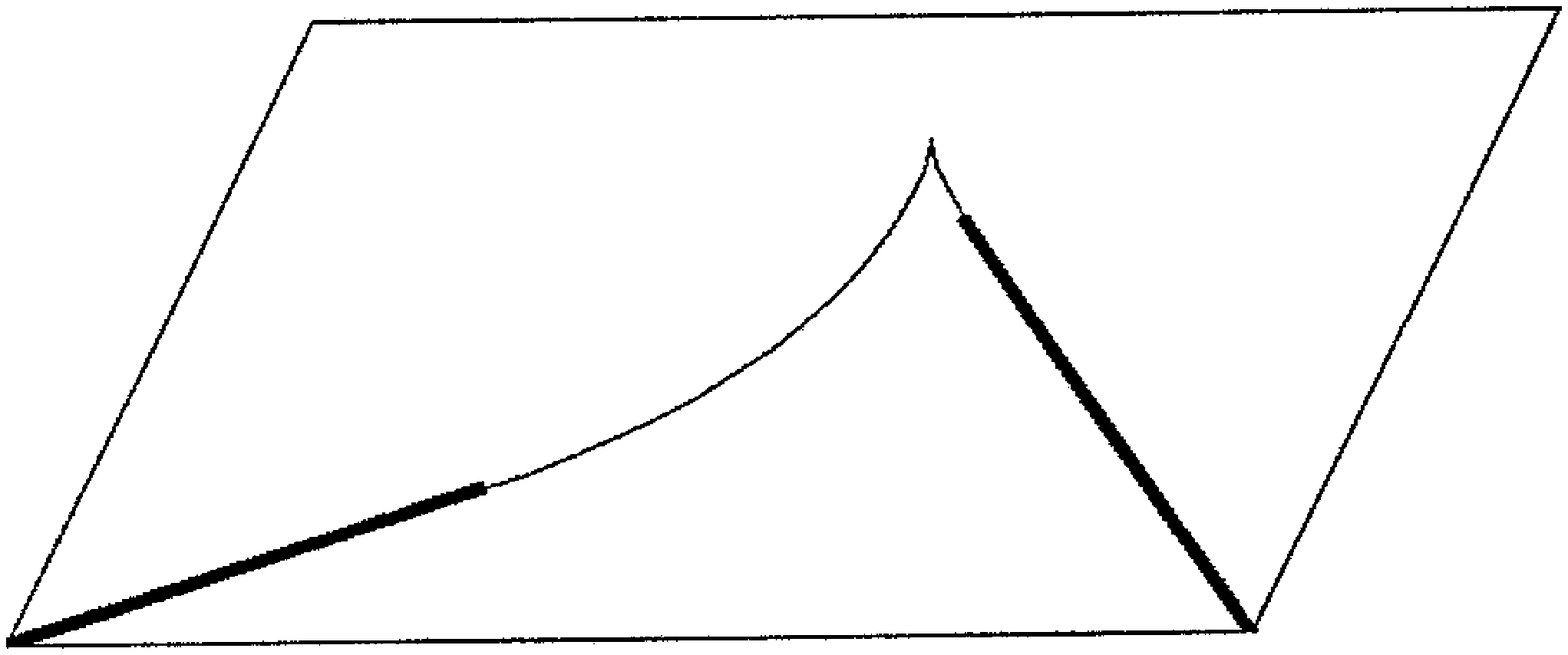}}
\centerline{Figure 11.}
\bigskip

\subheading{Final Observations}
Note that the condition which implies that the holonomy has an attracting fixed point 
is $\Im^2(\conj bd) > 1$. As the proof shows, this happens for parallelograms with 
area exactly $\pi\l^2$, and even for some with area slightly less. Under the 
assumption that $A \ge \pi\l^2$, the minimum value of 
$\,\Im(\conj bd) = \sinh|\b|\, \sinh|\d|\, \sin\theta\,$ is achieved for the square of
area $\pi\l^2$.

It seems intuitively clear that $z_+$ is the attracting fixed point and that $z_-$ is
the repelling fixed point.  To confirm this, let $h(z) = H{\cdot}z$ for
$z\in\C$. Since $h$ preserves $S^1$ and its orientation, $h'(z_+)$ and $h'(z_-)$
are both positive. A computation yields
$$
h'(z_\pm) = \( 2\Im(\conj bd)\(\Im(\conj bd) \pm \sqrt{\Im^2(\conj bd)-1}\) - 1 \)^{-2}.
$$
Remembering that $\Im(\conj bd) > 1$, it follows that $h'(z_+) < 1$ and $h'(z_-) > 1$,
justifying intuition.

As noted earlier, the results of Prytz and Hill imply that when the boundary of a 
small region is traced, the resulting holonomy has no fixed points.  There are,
however, short curves which result in holonomy with fixed points.  As an example,
let $v$ and $w$ be as in the proof (but without any assumption about the area of the
parallelogram), and consider the path consisting of the
parallelogram in the proof and its reflection through the origin. More specifically,
the path is the polygonal ``figure eight'' formed by the successive vertices: 0, $v$, 
$v+w$, $w$, $-w$, $-v-w$, $-v$, 0. Note that the oriented area bounded by this path is
0. The resulting holonomy is $\tilde H = e^X e^Y e^{-X}e^{-2Y}e^{-X}e^Y e^X$, where
$X$ and $Y$ are given by \thetag9. A computation shows that
$\tr\tilde H = 2 + 16 \Im^2(\conj bd) |ad + bc|^2$.
As long as $v$ and $w$ are independent (which implies $\Im(\conj bd) \ne 0$ and 
$ad + bc \ne 0$) we have $\tr\tilde H > 2$, and so $\tilde H$ has two fixed points.
Furthermore, the fixed points are antipodal. This isn't surprising, given the
symmetry of the ``figure eight'' about the base point 0. (The expressions for the
fixed points are very long and not very enlightening, so they are omitted.) 
Thus $\tilde H = e^Z$, where $Z = \o(u)$ for some $u\in\R^2$. As we have seen,
elements of $\su$ of this form generate $\su$.  It follows that the holonomy group
is generated by short loops.

%\newpage

\heading Figure Credits \endheading

Figures 1(a,b) are from ``Mathematical Machines'' by Francis J. Murray 
\cite{Mu, p.~348}, Copyright \copyright\ 1961 by Columbia University Press, and are 
reprinted with permission of the publisher. 
Figure~2 is from \cite{Pou}, with letters added.
Figure~7 is from \cite{G}.
Animated versions of some of the figures are available on the author's web page
(URL below).

%%% TXP %%% END Prytz56.tex %%%%%%%%%%%%%%%%%%%%%%%%%%%%%%%%%%%%%%%%%%%%%%%

%%% TXP %%% BEGIN PrytzRef.tex %%%%%%%%%%%%%%%%%%%%%%%%%%%%%%%%%%%%%%%%%%%%

\enddocument

     Back up possibilities for Fig 1.
	polar: Horsburgh, pg 207; Larson, pg 938; Granville, pg 448
        linear: possibly one of the detailed drawings in Horsburgh
                Meyer Zur Capellen, pg 183
        Schematic drawings of both in Lopshits, but this is copyright 1963 in USA
          You could easily make similar drawings.
%%% TXP %%% END PrytzRef.tex %%%%%%%%%%%%%%%%%%%%%%%%%%%%%%%%%%%%%%%%%%%%%%


\Refs

\widestnumber\key{KMS}

\ref\key Ba
\by G. Barnes
\paper Hatchet or Hacksaw Blade Planimeter
\jour Am. J. Physics
\vol 25 
\yr 1957
\pages 25--29
\endref

\ref\key Bo
\by W. M. Boothby
\book An Introduction to Differentiable Manifolds and Riemannian Geometry
\yr 1975
\publ Academic Press
\publaddr New York
\endref

\ref\key Cou
\by R. Courant
\book Differential and Integral Calculus
\bookinfo Vol. II
\yr 1934
\publ Nordemann
\publaddr New York
\endref

\ref\key Cox
\by H. S. M. Coxeter
\book Introduction to Geometry, 
\bookinfo second ed.
\yr 1989
\publ Wiley
\publaddr New York
\endref

\ref\key Cr
\by A. R. Crathorne 
\paper The Prytz Planimeter
\jour Am. Math. Monthly
\vol 15
\yr 1908
\pages 55--57
\endref

\ref\key Fe
\by M. Fecko
\paper Gauge-potential Approach to the Kinematics of a Moving Car
\jour Il Nuovo Cimento B
\vol 111 
\pages 1315--1332 
\yr 1996
\endref

\ref\key Fo1
\by R. L. Foote
\paper A Plenitude of Planimeters
\jour in preparation
\endref

\ref\key Fo2
\by R. L. Foote
\paper Planimeters and Isoperimetric Inequalities on Constant Curvature Surfaces
\jour in preparation
\endref

\ref\key G
\by J. Goodman (pub. anon.)
\paper Goodman's Hatchet Planimeter
\jour Engineering, Aug. 21, 1896
\pages 255--56
\endref

\ref\key He
\by O. Henrici
\paper Report on Planimeters
\jour British Assoc. for the Advancement of Science, Report of the 64th meeting
\yr 1894  
\pages 496--523
\endref

\ref\key Hi
\by F. W. Hill 
\paper The Hatchet Planimeter
\jour Philosophical Magazine, S.~5, Vol.~38, No.~ 232, Sept., 1894
\pages 265--269
\moreref
\jour Proc. of the Physical Society
\vol 13
\pages 229--234 \nofrills
\finalinfo (same paper appears twice)
\endref


\ref\key K
\by A. Kriloff
\paper On the Hatchet Planimeter
\jour Bulletin de l'Acad\'emie Imp\'eriale des Sciences de St.~P\'eters\-bourg, 
      T. XIX, No.~4 \& 5, Nov/Dec, 1903
\pages 221--227
\endref

\ref\key KN
\by S. Kobayashi and K. Nomizu
\book Foundations of Differential Geometry
\bookinfo Vol. I
\publ Wiley-Interscience \publaddr New York
\yr 1963
\endref

\ref\key KMS
\by I. Kol\'a\v r, P. W. Michor, and J. Slov\'ak
\book Natural Operations in Differential Geometry
\publ Springer-Verlag \publaddr Berlin
\yr 1993
\endref

\ref\key L
\by D. N. Lehmer
\paper Concerning the Tractrix of a Curve, with Planimetric Application
\jour Annals of Math
\vol 13 
\yr 1899
\pages 14--20
\endref

\ref\key Me
\by A. L. Menzin
\paper The Tractigraph, an Improved form of Hatchet Planimeter
\jour Engineering News,
Vol.~56, No.~6
\yr 1906
\pages 131--132
\endref

\ref\key Mo
\by F. Morley
\paper The ``No-Rolling'' Curves of Amsler's Planimeter
\jour Annals of Math
\vol 13 
\yr 1899
\pages 21--30
\endref

\ref\key Mu
\by F. J. Murray
\book Mathematical Machines, {\rm Vol.~2}, Analog Devices
\publ Columbia University Press
\publaddr New York
\yr 1961
\endref

\ref\key Pe
\by Olaf Pedersen 
\paper The Prytz Planimeter
\inbook From Ancient Omens to Statistical Mechanics
\eds J.L. Berggren and B.R. Goldstein
\publ University Library
\publaddr Copenhagen 
\yr 1987
\endref

\ref\key Poo
\by W. A. Poor
\book Differential Geometric Structures
\publ McGraw-Hill 
\yr 1981
\endref


\ref\key Pou
\by A.\ Poulain
\paper Les Aires des Tractrices et le Stang-Planim\`etre
\lang French
\jour J. de Math\'ematiques Sp\'eciales, Vol.~4, No.~2
\yr 1895
\pages 49--54
\endref

\ref\key Pr1
\by H. Prytz (pseud. `Z')
\paper Stangplanimetret
\jour Den Tekniske Forenings Tidsskrift
\lang Danish
\vol 10
\yr 1886
\pages 23--28 \nofrills
\finalinfo (appendix to Heinrich Ohrt, {\sl Om Planimetre,} 14--28)
\endref

\ref\key Pr2
\by H. Prytz
\paper The Hatchet Planimeter
\jour (letter to the editor), Engineering
\vol 57
\yr June 22, 1894
\page 813
\endref

\ref\key Pr3
\by H. Prytz
\paper The Prytz Planimeter 
\jour (two letters to the editor), Engineering
%\vol ?
\yr September 11, 1896
\page 347
\endref

\ref\key Pr4
\by H. Prytz
\paper The Hatchet Planimeter and `Tractigraph'
\jour (letter to the editor), Engineering News, Vol.~57, No.~14
\yr 1907
\page 386
\endref

\ref\key Sa
\by J. Satterly
\paper The Hatchet Planimeter
\jour J. Royal Astronomical Soc Canada, Vol.~15, No.~6
\yr 1921
\pages 221--243
\endref

\ref\key Sc
\by E. K. Scott
\paper An Improved Stang Planimeter
\jour Engineering
\yr Aug. 14, 1896
\pages 205--206
\endref   

\ref\key Se
\by C. L. Seigel
\book Topics in Complex Function Theory
\bookinfo Vol. II
\publ Wiley-Interscience \publaddr New York
\yr 1971
\endref   

\ref\key St
\by C. L. Strong
\paper An Excursion into the Problem of Measuring Irregular Areas
\jour Scientific American ({\sl The Amateur Scientist\/} column), Vol.~199, No.~2
\yr 1958
\pages 107--114 \nofrills
\finalinfo (letter from F. W. Niedenfuhr)
\endref


\endRefs
\enddocument